\documentclass[12pt]{article}
\usepackage{amsthm}
\usepackage{amssymb}
\usepackage{amsmath}
\usepackage{amsfonts}
\usepackage{enumerate}
\usepackage{verbatim}
\usepackage[dvips]{graphicx}
\usepackage[english]{babel}
\usepackage[cp1250]{inputenc}
\usepackage[T1]{fontenc}
\usepackage[a4paper,left=3cm,right=3cm,top=3cm,bottom=3cm]{geometry}

\newtheorem{thm}{Theorem}[section]
\newtheorem{lem}{Lemma}[section]

\newtheorem{cor}{Corollary}[section]
\newtheorem{df}{Definition}[section]

\newcommand{\Proof}{\noindent\textbf{{Proof:}}\newline}
\numberwithin{equation}{section}

\newcommand{\cbdu}{\quad\hfill\mbox{$\Box$}\\[3mm]}

\begin{document}

\title{\textbf{On Malliavin's differentiability of BSDE with time delayed generators driven by
Brownian motions and Poisson random measures}}

\author{\textbf{{\L}ukasz Delong}$^{1}$\footnote{corresponding author, tel/fax:(+48) 22 5648617.} , \textbf{Peter Imkeller}$^{2}$\\
\
\\
\footnotesize{Institute of Econometrics, Division of Probabilistic Methods}\\
\footnotesize{Warsaw School of Economics}\\
\footnotesize{Al. Niepodleglosci 162, 02-554 Warsaw, Poland}\\
\footnotesize{lukasz.delong@sgh.waw.pl}
\
\\
\
\\
\footnotesize{Institut f\"ur Mathematik}\\
\footnotesize{Humboldt-Universit\"at zu Berlin}\\
\footnotesize{Unter den Linden 6, 10099 Berlin, Germany}\\
\footnotesize{imkeller@mathematik.hu-berlin.de}}

\date{}

\maketitle

\newpage

\begin{abstract}
\noindent We investigate solutions of backward stochastic
differential equations (BSDE) with time delayed generators driven by
Brownian motions and Poisson random measures, that constitute the
two components of a L\'{e}vy process. In this new type of equations,
the generator can depend on the past values of a solution, by
feeding them back into the dynamics with a time lag. For such time
delayed BSDE, we prove existence and uniqueness of solutions
provided we restrict on a sufficiently small time horizon or the
generator possesses a sufficiently small Lipschitz constant. We
study differentiability in the variational or Malliavin sense and
derive equations that are satisfied by the Malliavin gradient
processes. On the chosen stochastic basis this addresses smoothness
both with respect to the continuous part of our L\'{e}vy process in
terms of the classical Malliavin derivative for Hilbert space valued
random variables, as well as with respect to the pure jump component
for which it takes the form of an increment quotient operator
related to the Picard difference operator.


\end{abstract}

\newpage

\section{Introduction}

\indent Introduced in \cite{PP}, backward stochastic differential
equations have been thoroughly studied in the literature during the
last decade, see \cite{K} or \cite{I1} and references therein.
Viewed from the perspective of Peng who interprets their key
structural feature as a \emph{nonlinear conditional expectation},
the close link to the stochastic calculus of variations or
Malliavin's calculus becomes apparent. In fact, in a Clark-Ocone
type formula, the control component of the solution pair of a BSDE
with a classical globally Lipschitz generator without time delay on
a Gaussian basis turns out to be the Malliavin trace of the other
component, see Proposition 5.3 in \cite{K} or Theorem 3.3.1 in
\cite{I1}. Not only this observation attributes an important role to
Malliavin's calculus in the context of stochastic control theory and
BSDE. As the simplest example, let us recall that hedging strategies
in complete market models corresponds to Malliavin derivatives of
wealth processes, see \cite{Ko}. The fine structure and sensitivity
properties of solutions of BSDE or systems of forward and backward
stochastic differential equations have been approached by means of
the stochastic calculus of variations (see \cite{AIdR1} and
\cite{AI}), and applied to provide explicit descriptions of delta
hedges of insurance related financial derivatives in \cite{AIdR2}.
Let us mention that Malliavin's calculus has been applied to prove
regularity of trajectories and thus to provide a first numerical
scheme for BSDE with generators of quadratic growth, see for
instance \cite{I2}. More generally, it has been established as a key
tool in the numerics of control theory and mathematical finance, for
instance to enhance the convergence speed of discretization schemes
for solutions of BSDE, see \cite {KK}, \cite{MT}. BSDE have proved
to be an efficient and powerful tool in a variety of applications in
stochastic control and mathematical finance. In all of these
applications, variational smoothness of their solutions is
fundamental for describing their properties.
\\
\indent In this spirit, and with the aim of clarifying smoothness in
the sense of the stochastic calculus of variations and related
properties of BSDE in a more general setting, in this paper we study
the equations with dynamics given for $t\in[0,T]$ by
$$Y(t)=\xi+\int_{t}^{T}f(s,Y_{s},Z_{s},U_{s})ds-\int_{t}^{T}Z(s)dW(s)-\int_{t}^{T}U(s,z)\tilde{M}(ds,dz).$$
An equation of this type will be called BSDE with time delayed
generator. It is driven by a L\'{e}vy process, the components of
which are given by a Brownian motion and a Poisson random measure.
In this new type of equations, a generator $f$ at time $s$ depends
in some measurable way on the past values of a solution
$(Y_{s},Z_{s},U_{s})=(Y(s+u),Z(s+u), U(s+u,.))_{-T\leq u\leq 0}$.
Very recently, time delayed BSDE driven by Brownian motion and with
Lipschitz continuous generators have been investigated for the the
first time in \cite{BI}, and in more depth in \cite{DI}. We would
like to refer the interested reader to the accompanying paper
\cite{DI}, where existence and uniqueness questions are treated, and
examples given in which multiple solutions or no solutions at all
exist. Further, several solution properties are investigated,
including the comparison principle, measure solutions, the
inheritance property of boundedness from terminal condition to
solution, as well as the \emph{BMO} martingale property for the
control component. We would like to point out that all results from
\cite{DI} can be extended and proved in the setting of
this paper.\\
\indent Our main findings are the following. First, we prove that a
unique solution exists, provided that the Lipschitz constant of the
generator is sufficiently small, or the equation is considered on a
sufficiently small time horizon. This is the extension of Theorem
2.1 from \cite{DI} to be expected. Secondly, we establish
Malliavin's differentiability of the solution of a time delayed
BSDE, both with respect to the continuous component of the L\'{e}vy
process, which coincides with the classical Malliavin derivative for
Hilbert-valued random variables, as well as with respect to the pure
jump part, in terms of an increment quotient operator related to
Picard's difference operator. We prove that the well-known
connection between $(Z,U)$ and the Malliavin trace of $Y$ still
holds in the case of time delayed generators. \\
\indent BSDE without time delays and driven by Poisson random
measures have already been thoroughly investigated in the
literature, see \cite{BBP}, \cite{Bech} or \cite{R}. But contrary to
the case with a Gaussian basis, smoothness results in the sense of
Malliavin's calculus have not been established yet in a systematic
way. To the best of our knowledge, only in \cite{Bo}, variational
differentiability of a solution of a forward-backward SDE with jumps
with respect to the Brownian component is considered while
differentiability with respect to the jump component is neglected.
\\
\indent We would like to emphasize that backward stochastic
differential equations with time-delayed generators arise in
financial and insurance problems dealing with pricing, hedging, risk
management and optimal control, see the working paper \cite{LD}. For
instance, they are encountered in the context of the optimal
liquidation problem of large trader's positions. A related optimal
control problem in terms of BSDE exhibits generators in which the
delayed feedback of the large trader's actions on the price dynamics
take the form of a delay effect in the sense considered in this
paper. As explained at the beginning of this section, Malliavin's
calculus plays a fundamental role in mathematical finance and
optimal control. We believe that the results concerning Malliavin's
differentiability obtained in this paper are as important to
describe parameter sensitivity properties of financial derivatives
in this more general setting as they are in \cite{AIdR2} for
generalizing the Black-Scholes delta hedge to incomplete markets in
a purely probabilistic approach via BSDE. \\
\indent This paper is structured as follows. Section 2 deals with
the existence and uniqueness problem. In Section 3 we survey
concepts of the canonical L\'{e}vy space and variational
differentiation, and prove some technical lemmas. The main theorem
concerning Malliavin smoothness of a solution, and the
interpretation of the latter in terms of a Malliavin trace is proved
in Section 4.

\section{Existence and uniqueness of a solution}
\indent We consider a probability space
$(\Omega,\mathcal{F},\mathbb{P})$ with a filtration
$\mathbb{F}=(\mathcal{F}_{t})_{0\leq t\leq T}$, where $T<\infty$ is
a finite time horizon. We assume that the filtration $\mathbb{F}$ is
the natural filtration generated by a L\'{e}vy process
$L:=(L(t),0\leq t\leq T)$ and that $\mathcal{F}_{0}$ contains all
sets of $\mathbb{P}$-measure zero, so that the \emph{usual
conditions} are fulfilled. As usual, by $\mathcal{B}(X)$ we denote
the Borel sets of a topological space $X$, while $\lambda$ stands for Lebesgue measure.\\
\indent It is well-known that a L\'{e}vy process satisfies the
L\'{e}vy-It\^{o} decomposition
\begin{eqnarray*}
L(t)=at+\sigma W(t)+\int_{0}^{t}\int_{|z|\geq 1}z N(ds,
dz)+\int_{0}^{t}\int_{0<|z|< 1}z (N(ds, dz)-\nu(dz)ds),
\end{eqnarray*}
for $0\leq t\leq T$, with $a\in\mathbb{R},\sigma\geq 0$. Here
$W:=(W(t),0\leq t\leq T)$ denotes a Brownian motion and $N$ a random
measure on $[0,T]\times(\mathbb{R}-\{0\})$, so that $W$ and $N$ are
independent. The random measure $N$
\begin{eqnarray*}
N(t,A)=\sharp\{0\leq s\leq t;\triangle L(s)\in A\},\quad 0\leq t\leq
T, A\in\mathcal{B}(\mathbb{R}-\{0\}),
\end{eqnarray*}
counts the number of jumps of a given size. It is called Poisson
random measure since, for $t\in[0,T]$ and a Borel set $A$ such that
its closure does not contain zero, $N(t,A)$ is a Poisson distributed
random variable. The $\sigma$-finite measure $\nu$, defined on
$\mathcal{B}(\mathbb{R}-\{0\})$, appears in the compensator
$\lambda\otimes \nu$ of the random measure $N$. The compensated
Poisson random measure (or martingale-valued measure) is denoted by
$\tilde{N}(t,A)=N(t,A)-t\nu(A)$,
$t\in[0,T],A\in\mathcal{B}(\mathbb{R}-\{0\})$. In this paper we deal
with the random measure
\begin{eqnarray*}
\tilde{M}(t,A)&=&\int_{0}^{t}\int_{A}z\tilde{N}(ds,dz)\nonumber\\
&=&
\int_{0}^{t}\int_{A}zN(ds,dz)-\int_{0}^{t}\int_{A}z\nu(dz)ds,\quad
0\leq t\leq T, A\in\mathcal{B}(\mathbb{R}-\{0\}).
\end{eqnarray*}
It can be considered as a compensated compound Poisson random
measure as, for a fixed $t\in[0,T]$ and a Borel set $A$ the closure
of which does not contain zero, $\int_{0}^{t}\int_{A}zN(ds,dz)$ is a
compound Poisson distributed random variable. Finally, we introduce
the $\sigma$-finite measure
\begin{eqnarray*}
m(A)=\int_{A}z^{2}\nu(dz),\quad A\in\mathcal{B}(\mathbb{R}-\{0\}).
\end{eqnarray*}
\noindent For details concerning L\'{e}vy processes, Poisson random
measures and integration with respect to martingale-valued random
measures we refer
the reader to Chapter 2 and Chapter 4 of \cite{A}.\\
\indent Let us now turn to the main subject of this paper. We study
solutions $(Y,Z,U):=(Y(t),Z(t),U(t,z))_{0\leq t\leq
T,z\in(\mathbb{R}-\{0\})}$ of a BSDE with time delayed generator,
the dynamics of which is given by
\begin{eqnarray}\label{bsde}
Y(t)&=&\xi+\int_{t}^{T}f(s,Y_{s},Z_{s},U_{s})ds\nonumber\\
&&-\int_{t}^{T}Z(s)dW(s)-\int_{t}^{T}\int_{\mathbb{R}-\{0\}}U(s,z)\tilde{M}(ds,dz),\quad
0\leq t\leq T.
\end{eqnarray}
The generator $f$ depends on the past values of the solution, fed
back into the system with a time delay, denoted by
$Y_{s}:=(Y(s+v))_{-T\leq v\leq 0}, Z_{s}:=(Z(s+v))_{-T\leq v\leq 0}$
and $U_{s}:=(U(s+v,.))_{-T\leq v\leq 0}, 0\le s\le T$. We always set
$Z(t)=U(t,.)=0$ and $Y(t)=Y(0)$ for $t<0$. Note that the measure
$\tilde{M}$, not $\tilde{N}$, is taken to drive the jump noise. The
reason for this is that we adopt the concepts of Malliavin calculus
on the canonical L\'{e}vy space from \cite{Sole}, which is
formulated in terms of multiple stochastic integrals with respect to
$\tilde{M}$.
\\
\indent We shall work with the function spaces of the following
definition.
\begin{df}
\begin{itemize}
    \item[1.] Let $L^{2}_{-T}(\mathbb{R})$ denote the space of measurable
        functions $z:[-T,0]\rightarrow\mathbb{R}$ satisfying
        $$\int_{-T}^{0}|z(t)|^{2}dt<\infty.$$
        \item[2.] Let $L^{2}_{-T,m}(\mathbb{R})$ denote the space of product measurable
        functions $u:[-T,0]\times(\mathbb{R}-\{0\})\rightarrow\mathbb{R}$ satisfying
        $$\int_{-T}^{0}\int_{\mathbb{R}-\{0\}}|u(t,z)|^{2}m(dz)dt<\infty.$$
    \item[3.] Let $L^{\infty}_{-T}(\mathbb{R})$ denote the space of bounded, measurable functions
        $y:[-T,0]\rightarrow\mathbb{R}$ such that
        $$\sup_{t\in[-T,0]} \big|y(t)\big|^{2}<\infty.$$
    \item[4.] Let $\mathbb{L}^{2}(\mathbb{R)}$ denote
        the space of $\mathcal{F}_{T}$-measurable random variables $\xi:\Omega\rightarrow\mathbb{R}$
        which fulfill
        $$\mathbb{E}\big[\big|\xi\big|^{2}\big]<\infty.$$
    \item [5.]Let $\mathbb{H}^{2}_{T}(\mathbb{R)}$ denote
        the space of predictable processes $Z:\Omega\times[0,T]\rightarrow\mathbb{R}$
        such that
        $$\mathbb{E}\big[\int_{0}^{T}\big|Z(t)\big|^{2}dt\big]<\infty.$$
    \item [6.]Let $\mathbb{H}^{2}_{T,m}(\mathbb{R)}$ denote
        the space of predictable processes $U:\Omega\times[0,T]\times(\mathbb{R}-\{0\})\rightarrow\mathbb{R}$ satisfying
        $$\mathbb{E}\big[\int_{0}^{T}\int_{\mathbb{R}-\{0\}}\big|U(t,z)\big|^{2}m(dz)dt\big]<\infty.$$
    \item [7.] Finally, let $\mathbb{S}^{2}_{T}(\mathbb{R)}$ denote
        the space of $\mathbb{F}$-adapted, product measurable processes
        $Y:\Omega\times[0,T]\rightarrow\mathbb{R}$ satisfying
        $$\mathbb{E}\big[\sup_{t\in[0,T]} \big|Y(t)\big|^{2}\big]<\infty.$$
  \end{itemize}
The spaces
$\mathbb{H}^{2}_{T}(\mathbb{R)}$,$\mathbb{H}^{2}_{T,m}(\mathbb{R)}$
and $\mathbb{S}^{2}_{T}(\mathbb{R)}$ are endowed with the norms
\begin{eqnarray*}
\big\|Z\big\|^{2}_{\mathbb{H}^{2}_{T}}&=&\mathbb{E}\big[\int_{0}^{T}e^{\beta
t}\big|Z(t)\big|^{2}dt\big],\\
\big\|U\big\|^{2}_{\mathbb{H}^{2}_{T,m}}&=&\mathbb{E}\big[\int_{0}^{T}\int_{\mathbb{R}-\{0\}}e^{\beta
t}\big|U(t,z)\big|^{2}m(dz)dt\big],\\
\big\|Y\big\|^{2}_{\mathbb{S}^{2}_{T}}&=&\mathbb{E}\big[\sup_{t\in[0,T]}
e^{\beta t}\big|Y(t)\big|^{2}\big],
\end{eqnarray*}
with some $\beta>0$.
\end{df}
\noindent Predictability of $Z$ means measurability with respect to
the predictable $\sigma$-algebra, which we denote by $\mathcal{P}$,
while predictability of $U$ means measurability with respect to the
product $\mathcal{P}\otimes\mathcal{B}(\mathbb{R}-\{0\})$. In the
sequel let us simply write
$\mathbb{S}^{2}(\mathbb{R})\times\mathbb{H}^{2}(\mathbb{R})\times\mathbb{H}^{2}_{m}(\mathbb{R})$
instead of
$\mathbb{S}^{2}_{T}(\mathbb{R})\times\mathbb{H}^{2}_{T}(\mathbb{R})\times\mathbb{H}^{2}_{T,m}(\mathbb{R})$.\\
\indent We start with establishing existence and uniqueness of a
solution of \eqref{bsde} under the following hypotheses:
\begin{description}
  \item[(A1)] the terminal value $\xi\in\mathbb{L}^{2}(\mathbb{R})$,
  \item[(A2)] $m$ is a finite measure, i.e. $\int_{\mathbb{R}-\{0\}}z^{2}\nu(dz)<\infty$,
  \item[(A3)] the generator
  $f:\Omega\times [0,T]\times L^{\infty}_{-T}(\mathbb{R}) \times L^{2}_{-T}(\mathbb{R})\times L^{2}_{-T,m}(\mathbb{R})\rightarrow \mathbb{R}$
  is product measurable, $\mathbb{F}$-adapted and Lipschitz
  continuous in the sense that for a probability measure $\alpha$ on
  $([-T,0],\mathcal{B}([-T,0]))$ and with a constant $K>0$
  \begin{eqnarray*}
  \lefteqn{|f(\omega,t,y_{t},z_{t},u_{t})-f(\omega,t,\tilde{y}_{t},\tilde{z}_{t},\tilde{u}_{t})|^{2}}\nonumber\\
  &\leq&
  K\Big(\int_{-T}^{0}|y(t+v)-\tilde{y}(t+v)|^{2}\alpha(dv)+\int_{-T}^{0}|z(t+v)-\tilde{z}(t+v)|^{2}\alpha(dv)\nonumber\\
  &&+\int_{-T}^{0}\int_{\mathbb{R}-\{0\}}|u(t+v,z)-\tilde{u}(t+v,z)|^{2}m(dz)\alpha(dv)\Big),
  \end{eqnarray*}
  holds for $\mathbb{P}\otimes\lambda$-a.e. $(\omega,t)\in\Omega\times[0,T]$, for any
  $(y_{t},z_{t},u_{t}),(\tilde{y}_{t},\tilde{z}_{t},\tilde{u}_{t})\in L^{\infty}_{-T}(\mathbb{R})
  \times L^{2}_{-T}(\mathbb{R})\times L^{2}_{-T,m}(\mathbb{R})$,
  \item[(A4)] $\mathbb{E}\big[\int_{0}^{T}|f(t,0,0,0)|^{2}dt\big]<\infty$,
  \item[(A5)] $f(\omega, t,.,.,.)=0$ for $\omega\in\Omega, t<0$.
\end{description}
For convenience, in the notation of $f$ the dependence on $\omega$
is omitted and we write $f(t,.,.,.)$ for $f(\omega,t,.,.,.)$ etc. We
remark that $f(t,0,0,0)$ in \textbf{(A4)} should be understood as
the value of the generator $f(t,y_{t},z_{t},u_{t})$ at
$y_t=z_t=u_t=0$. We would like to point out that assumption
\textbf{(A5)} in fact allows us to take $Y(t)=Y(0)$ and
$Z(t)=U(t,.)=0$ for $t<0$ as a solution of \eqref{bsde}. Finally,
let us recall that under \textbf{(A2)} and for an integrand
$U\in\mathbb{H}_{m}^{2}(\mathbb{R})$, the stochastic integral with
respect to the martingale-valued measure $\tilde{M}$
\begin{eqnarray*}
\int_{0}^{t}\int_{\mathbb{R}-\{0\}}U(s,z)\tilde{M}(ds,dz),\quad
0\leq t\leq T,
\end{eqnarray*}
is well-defined in the It\^{o} sense, see Chapter 4.1 in \cite{A}.\\
\indent First let us notice that for
$(Y,Z,U)\in\mathbb{S}^{2}(\mathbb{R})\times\mathbb{H}^{2}(\mathbb{R})\times\mathbb{H}^{2}_{m}(\mathbb{R})$
the generator is well-defined and integrable as a consequence of
\begin{eqnarray}\label{estf}
\lefteqn{\int_{0}^{T}|f(t,Y_{t},Z_{t},U_{t})|^{2}dt\leq2\int_{0}^{T}|f(t,0,0,0)|^{2}dt
+2K(\int_{0}^{T}\int_{-T}^{0}|Y(t+v)|^{2}\alpha(dv)dt}\nonumber\\
&&+\int_{0}^{T}\int_{-T}^{0}|Z(t+v)|^{2}\alpha(dv)dt
+\int_{0}^{T}\int_{-T}^{0}\int_{\mathbb{R}-\{0\}}|U(t+v,z)|^{2}m(dz)\alpha(dv)dt)\nonumber\\
&=&2\int_{0}^{T}|f(t,0,0,0)|^{2}dt+2K\int_{-T}^{0}\int_{v}^{T+v}|Y(w)|^{2}dw\alpha(dv)\nonumber\\
&&+2K\int_{-T}^{0}\int_{v}^{T+v}|Z(w)|^{2}dw\alpha(dv)\nonumber\\
&&+2K\int_{-T}^{0}\int_{v}^{T+v}\int_{\mathbb{R}-\{0\}}|U(w,z)|^{2}m(dz)dw\alpha(dv)\nonumber\\
&\leq& 2\int_{0}^{T}|f(t,0,0,0)|^{2}dt+2K\big(T\sup_{w\in[0,T]}|Y(w)|^{2}\nonumber\\
&&+\int_{0}^{T}|Z(w)|^{2}dw+\int_{0}^{T}\int_{\mathbb{R}-\{0\}}|U(w,z)|^{2}m(dz)dw\big)<\infty,\quad
\mathbb{P}-a.s.,
\end{eqnarray}
where we apply \textbf{(A3)}, Fubini's theorem, use the assumption
that $Z(t)=U(t,.)=0$
and $Y(t)=Y(0)$ for $t<0$ and the fact that the measure $\alpha$ is a probability measure.\\
\indent The main theorem of this section is an extension of Theorem
2.1 from \cite{DI}. Although the extension is quite natural, the
proof is given for completeness and convenience of the
reader. The key result follows from the following a priori estimates. \\
\begin{lem}\label{priori}
Let $(Y,Z,U), (\tilde{Y},\tilde{Z},\tilde{U})
\in\mathbb{S}^{2}(\mathbb{R})\times\mathbb{H}^{2}(\mathbb{R})\times\mathbb{H}^{2}_{m}(\mathbb{R})$
denote solutions of \eqref{bsde} with corresponding parameters
$(\xi,f)$ and $(\tilde{\xi},\tilde{f})$ which satisfy the
assumptions \textbf{(A1)}-\textbf{(A5)}. Then the following
inequalities hold
\begin{eqnarray}\label{p1}
\lefteqn{\|Z-\tilde{Z}\|^{2}_{\mathbb{H}^{2}}+\|U-\tilde{U}\|^{2}_{\mathbb{H}^{2}_{m}}}\nonumber\\
&\leq& e^{\beta T}\mathbb{E}\big[\big|\xi-\tilde{\xi}\big|^{2}\big]+
\frac{1}{\beta}\mathbb{E}\big[\int_{0}^{T}e^{\beta
t}|f(t,Y_{t},Z_{t},U_{t})-\tilde{f}(t,\tilde{Y}_{t},\tilde{Z}_{t},\tilde{U}_{t})|^{2}dt\big],
\end{eqnarray}
\begin{eqnarray}\label{p2}
\lefteqn{\|Y-\tilde{Y}\|^{2}_{\mathbb{S}^{2}}}\\\nonumber &\leq&
8e^{\beta T}\mathbb{E}\big[\big|\xi-\tilde{\xi}\big|^{2}\big] +
8T\mathbb{E}\big[\int_{0}^{T}e^{\beta
t}|f(t,Y_{t},Z_{t},U_{t})-\tilde{f}(t,\tilde{Y}_{t},\tilde{Z}_{t},\tilde{U}_{t})|^{2}dt\big].
\end{eqnarray}
\end{lem}
\Proof The inequality \eqref{p1} follows by a straightforward
extension of Lemma 3.2.1 from \cite{I1}, by only adding an
additional stochastic integral with respect to $\tilde{M}$. In order
to prove the second inequality, first notice that for $t\in[0,T]$
\begin{eqnarray*}
Y(t)-\tilde{Y}(t)=
\mathbb{E}\big[\xi-\tilde{\xi}+\int_{t}^{T}(f(s,Y_{s},Z_{s},U_{s})
-\tilde{f}(s,\tilde{Y}_{s},\tilde{Z}_{s},\tilde{U}_{s}))ds|\mathcal{F}_{t}\big],
\end{eqnarray*}
and
\begin{eqnarray*}
\lefteqn{e^{\frac{\beta}{2}t}|Y(t)-\tilde{Y}(t)|}\nonumber\\
&\leq&
e^{\frac{\beta}{2}T}\mathbb{E}\big[\big|\xi-\tilde{\xi}\big||\mathcal{F}_{t}\big]+
\mathbb{E}\big[\int_{0}^{T}e^{\frac{\beta}{2}
s}|f(s,Y_{s},Z_{s},U_{s})-\tilde{f}(s,\tilde{Y}_{s},\tilde{Z}_{s},\tilde{U}_{s})|ds|\mathcal{F}_{t}\big],
\end{eqnarray*}
hold $\mathbb{P}$-a.s.. Doob's martingale inequality and
Cauchy-Schwarz' inequality yield the second  estimate. The reader
may also consult Proposition 2.2 in \cite{BBP}
or Proposition 3.3 in \cite{Bech}, where similar estimates for BSDE with jumps are derived.\cbdu \\

\begin{thm}\label{thmain}
Assume that \textbf{(A1)}-\textbf{(A5)} hold. For a sufficiently
small time horizon $T$ or for a sufficiently small Lipschitz
constant $K$ of the generator $f$, more precisely if for some
$\beta>0$
\begin{eqnarray*}
\delta(T,K,\beta,\alpha):=(8T+\frac{1}{\beta})K\int_{-T}^{0}e^{-\beta
v}\alpha(dv)\max\{1,T\}<1,
\end{eqnarray*}
the backward stochastic differential equation \eqref{bsde} has a
unique solution $(Y,Z,U)\in\mathbb{S}^{2}(\mathbb{R})\times
\mathbb{H}^{2}(\mathbb{R})\times\mathbb{H}^{2}_{m}(\mathbb{R})$.
\end{thm}
\Proof We follow the classical Picard type iteration scheme (see
Theorem 2.1 in \cite{K} or Theorem 3.2.1 in \cite{I1}) to prove
existence and uniqueness of a solution.\\
\indent Let $Y^{0}(t)=Z^{0}(t)=U^{0}(t,z)=0$,
$(t,z)\in[0,T]\times(\mathbb{R}-\{0\})$.

\noindent Step 1) We show that the recursive definition
\begin{eqnarray}\label{bsden}
\lefteqn{Y^{n+1}(t)=\xi+\int_{t}^{T}f(s,Y^{n}_{s},Z^{n}_{s},U^{n}_{s})ds}\nonumber\\
&&-\int_{t}^{T}Z^{n+1}(s)dW(s)-\int_{t}^{T}\int_{\mathbb{R}-\{0\}}U^{n+1}(s,z)\tilde{M}(ds,dz),\quad
0\leq t\leq T,
\end{eqnarray}
makes sense. More precisely, we show that given
$(Y^{n},Z^{n},U^{n})\in\mathbb{S}^{2}(\mathbb{R})\times
\mathbb{H}^{2}(\mathbb{R})\times\mathbb{H}^{2}_{m}(\mathbb{R})$,
equation \eqref{bsden} has a unique solution
$(Y^{n+1},Z^{n+1},U^{n+1})\in\mathbb{S}^{2}(\mathbb{R})\times
\mathbb{H}^{2}(\mathbb{R})\times\mathbb{H}^{2}_{m}(\mathbb{R})$.\\
\noindent Applying inequality \eqref{estf}, we can conclude that
\begin{eqnarray*}
\mathbb{E}\big[\int_{0}^{T}|f(t,Y^{n}_{t},Z^{n}_{t},U^{n}_{t})|^{2}dt\big]&\leq&2\mathbb{E}\big[\int_{0}^{T}|f(t,0,0,0)|^{2}dt\big]\nonumber\\
&&+2K\big(T\|Y^{n}\|_{\mathbb{S}^{2}}+\|Z^{n}\|_{\mathbb{H}^{2}}+\|U^{n}\|_{\mathbb{H}^{2}_{m}}\big)<\infty.
\end{eqnarray*}
As in the case of BSDE without time delays, the martingale
representation, see Theorem 13.49 in \cite{Y}, provides a unique
process $Z^{n+1}\in\mathbb{H}^{2}(\mathbb{R})$ and a unique
predictable process $\bar{U}^{n+1}$ satisfying
\begin{eqnarray*}
\mathbb{E}\big[\int_{0}^{T}\int_{\mathbb{R}-\{0\}}|\bar{U}^{n+1}(t,z)|^{2}\nu(dz)dt\big]<\infty,
\end{eqnarray*}
so that
\begin{eqnarray*}
\lefteqn{\xi+\int_{0}^{T}f(t,Y^{n}_{t},Z^{n}_{t},U^{n}_{t})dt=\mathbb{E}\big[\xi+\int_{0}^{T}f(t,Y^{n}_{t},Z^{n}_{t},U^{n}_{t})dt\big]}\nonumber\\
&&+\int_{0}^{T}Z^{n+1}(t)dW(t)+\int_{0}^{T}\int_{\mathbb{R}-\{0\}}\bar{U}^{n+1}(t,z)\tilde{N}(dt,dz),\quad
\mathbb{P}-a.s.
\end{eqnarray*}
For $(t,z)\in[0,T]\times(\mathbb{R}-\{0\})$ we get
$U^{n+1}(t,z)=\frac{\bar{U}^{n+1}(t,z)}{z}\in\mathbb{H}^{2}_{m}(\mathbb{R}),$
and have the required representation
\begin{eqnarray*}
\lefteqn{\xi+\int_{0}^{T}f(t,Y^{n}_{t},Z^{n}_{t},U^{n}_{t})dt=\mathbb{E}\big[\xi+\int_{0}^{T}f(t,Y^{n}_{t},Z^{n}_{t},U^{n}_{t})dt\big]}\nonumber\\
&&+\int_{0}^{T}Z^{n+1}(t)dW(t)+\int_{0}^{T}\int_{\mathbb{R}-\{0\}}U^{n+1}(t,z)(t)\tilde{M}(dt,dz),\quad
\mathbb{P}-a.s..
\end{eqnarray*}
Finally, we take $Y^{n+1}$ as a progressively measurable,
c\`{a}dl\`{a}g modification of
\begin{eqnarray*}
Y^{n+1}(t)(\omega)=
\mathbb{E}\big[\xi+\int_{t}^{T}f(s,Y^{n}_{s},Z^{n}_{s},U^{n}_{s}))ds|\mathcal{F}_{t}\big],\quad
\omega \in\Omega, t\in[0,T].
\end{eqnarray*}
Similarly as in Lemma \ref{priori}, Doob's martingale inequality,
Cauchy-Schwarz' inequality and the estimates \eqref{estf} yield that
$Y^{n+1}\in\mathbb{S}^{2}(\mathbb{R})$.\\

\noindent Step 2) We prove the convergence of the sequence
$(Y^{n},Z^{n},U^{n})$ in
$\mathbb{S}^{2}(\mathbb{R})\times\mathbb{H}^{2}(\mathbb{R})\times\mathbb{H}^{2}_{m}(\mathbb{R})$.\\
\noindent The estimates \eqref{p1} and \eqref{p2} provide the
inequality
\begin{eqnarray}\label{contr1}
\lefteqn{\big\|Y^{n+1}-Y^{n}\big\|^{2}_{\mathbb{S}^{2}}+\big\|Z^{n+1}-Z^{n}\big\|^{2}_{\mathbb{H}^{2}}+\big\|U^{n+1}-U^{n}\big\|^{2}_{\mathbb{H}^{2}_{m}}}\nonumber\\
&&\leq(8T+\frac{1}{\beta})\mathbb{E}\big[\int_{0}^{T}e^{\beta
t}|f(t,Y^{n}_{t},Z^{n}_{t},U^{n}_{t})-f(t,Y^{n-1}_{t},Z^{n-1}_{t},U^{n-1}_{t})|^{2}dt\big].
\end{eqnarray}
By applying the Lipschitz condition \textbf{(A3)}, Fubini's theorem,
changing variables and using the assumption $\forall n\geq 0$
$Y^{n}(s)=Y^{n}(0)$ and $Z^{n}(s)=U^{n}(s,.)=0$ for $s<0$, we can
derive
\begin{eqnarray}\label{contr2}
\lefteqn{\mathbb{E}\big[\int_{0}^{T}e^{\beta t}|f(t,Y^{n}_{t},Z^{n}_{t},U^{n}_{t})-f(t,Y^{n-1}_{t},Z^{n-1}_{t},U^{n-1}_{t})|^{2}dt\big]}\nonumber\\
&\leq& K\mathbb{E}\big[\int_{0}^{T}e^{\beta
t}\int_{-T}^{0}|Y^n(t+v)-Y^{n-1}(t+v)|^{2}\alpha(dv)dt\nonumber\\
&&+\int_{0}^{T}e^{\beta t}\int_{-T}^{0}|Z^n(t+v)-Z^{n-1}(t+v)|^{2}\alpha(dv)dt\nonumber\\
&&+\int_{0}^{T}e^{\beta t}\int_{-T}^{0}\int_{\mathbb{R}-\{0\}}|U^n(t+v,z)-U^{n-1}(t+v,z)|^{2}m(dz)\alpha(dv)dt\big]\nonumber\\
&=&K\mathbb{E}\big[\int_{-T}^{0}e^{-\beta v}\int_{0}^{T}e^{\beta (t+v)}|Y^n(t+v)-Y^{n-1}(t+v)|^{2}dt\alpha(dv)\nonumber\\
&&+\int_{-T}^{0}e^{-\beta v}\int_{0}^{T}e^{\beta (t+v)}|Z^n(t+v)-Z^{n-1}(t+v)|^{2}dt\alpha(dv)\nonumber\\
&&+\int_{-T}^{0}e^{-\beta v}\int_{0}^{T}\int_{\mathbb{R}-\{0\}}e^{\beta (t+v)}|U^n(t+v,z)-U^{n-1}(t+v,z)|^{2}m(dz)dt\alpha(dv)\big]\nonumber\\
&=&K\mathbb{E}\big[\int_{-T}^{0}e^{-\beta v}\int_{v}^{T+v}e^{\beta w}|Y^n(w)-Y^{n-1}(w)|^{2}dw\alpha(dv)\nonumber\\
&&+\int_{-T}^{0}e^{-\beta v}\int_{v}^{T+v}e^{\beta w}|Z^n(w)-Z^{n-1}(w)|^{2}dw\alpha(dv)\nonumber\\
&&+\int_{-T}^{0}e^{-\beta v}\int_{v}^{T+v}\int_{\mathbb{R}-\{0\}}e^{\beta w}|U^n(w,z)-U^{n-1}(w,z)|^{2}m(dz)dw\alpha(dv)\big]\nonumber\\
&\leq&K\int_{-T}^{0}e^{-\beta
v}\alpha(dv)\big(T\big\|Y^{n}-Y^{n-1}\big\|^{2}_{\mathbb{S}^{2}}+\big\|Z^{n}-Z^{n-1}\big\|^{2}_{\mathbb{H}^{2}}+
\big\|U^{n}-U^{n-1}\big\|^{2}_{\mathbb{H}^{2}_{m}}\big).\nonumber\\
\end{eqnarray}
From \eqref{contr1} and \eqref{contr2}, we obtain
\begin{eqnarray}\label{contr}
\lefteqn{\big\|Y^{n+1}-Y^{n}\big\|^{2}_{\mathbb{S}^{2}}+\big\|Z^{n+1}-Z^{n}\big\|^{2}_{\mathbb{H}^{2}}+\big\|U^{n+1}-U^{n}\big\|^{2}_{\mathbb{H}^{2}_{m}}}\nonumber\\
&&\leq
\delta(T,K,\beta,\alpha)\big(\big\|Y^{n}-Y^{n-1}\big\|^{2}_{\mathbb{S}^{2}}+\big\|Z^{n}-Z^{n-1}\big\|^{2}_{\mathbb{H}^{2}}+\big\|U^{n}-U^{n-1}\big\|^{2}_{\mathbb{H}^{2}_{m}}),
\end{eqnarray}
with
$$\delta(T,K,\beta,\alpha)=(8T+\frac{1}{\beta})K\int_{-T}^{0}e^{-\beta v}\alpha(dv)\max\{1,T\}.$$
For $\beta=\frac{1}{T}$ we have
\begin{eqnarray*}
\delta(T,K,\beta,\alpha)\leq 9TKe\max\{1,T\}.
\end{eqnarray*}
For sufficiently small $T$ or sufficiently small $K$, the inequality
\eqref{contr} provides a unique limit
$(Y,Z,U)\in\mathbb{S}^{2}(\mathbb{R})\times\mathbb{H}^{2}(\mathbb{R})\times\mathbb{H}^{2}_{m}(\mathbb{R})$
of the converging sequence $(Y^{n},Z^{n},U^{n})_{n\in\mathbb{N}}$,
which satisfies the fixed point equation
\begin{eqnarray*}
Y(t)=
\mathbb{E}\big[\xi+\int_{t}^{T}f(s,Y_{s},Z_{s},U_{s})ds|\mathcal{F}_{t}\big],\quad
\mathbb{P}-a.s., 0\leq t\leq T.
\end{eqnarray*}
\noindent Step 3) We define the solution component $\bar{Y}$ of
\eqref{bsdem} as a progressively measurable, c\`{a}dl\`{a}g
modification of
\begin{eqnarray*}
\bar{Y}(t)(\omega)=
\mathbb{E}\big[\xi+\int_{t}^{T}f(s,Y_{s},Z_{s},U_{s})ds|\mathcal{F}_{t}\big],\quad
\omega\in\Omega, t\in[0,T],
\end{eqnarray*}
where $(Y,Z,U)$ is the limit constructed in Step 2). \cbdu

\indent We point out that in general under the assumptions
\textbf{(A1)}-\textbf{(A5)}, existence and uniqueness may fail to
hold for bigger time horizon $T$ or bigger Lipschitz constant $K$.
See \cite{DI} for examples. However, for some special classes of
generators existence and uniqueness may be proved for an arbitrary
time horizon and for arbitrary global Lipschitz constants. These
include generators independent of $y$ with a delay measure $\alpha$
supported on $[-\gamma,0]$ with a sufficiently small time delay
$\gamma$, following Theorem 2.2 in \cite{DI}, or generators
considered in \cite{BI} consisting of separate components in $z$ and
$u$, following Theorem 1 in \cite{BI}.

\section{Malliavin's calculus for canonical L\'{e}vy processes}
There are various ways to develop Malliavin's calculus for L\'{e}vy
processes. In this paper we adopt the approach from \cite{Sole}
based on a chaos decomposition in terms of multiple stochastic
integrals with respect to the random measure $\tilde{M}$. In this
setting, we will construct a suitable canonical space, on which a
variational derivative with respect to the pure jump part of a
L\'{e}vy
process can be computed in a pathwise sense.\\
\indent In this section we give an overview of the approach of
Malliavin's calculus on canonical L\'{e}vy space according to
\cite{Sole} (see \cite{Sole} for details). We then prove some
technical results concerning the commutation of integration and
variational differentiation,
which are needed in the next section.\\
\indent We assume that the probability space
$(\Omega,\mathcal{F},\mathbb{P})$ is the product of two canonical
spaces
$(\Omega_{W}\times\Omega_{N},\mathcal{F}_{W}\times\mathcal{F}_{N},\mathbb{P}_{W}\times\mathbb{P}_{N})$,
and the filtration $\mathbb{F} = (\mathcal{F}_t)_{t\in[0,T]}$ the
canonical filtration completed for $\mathbb{P}.$ The space
$(\Omega_{W},\mathcal{F}_{W},\mathbb{P}_{W})$ is the usual canonical
space for a one-dimensional Brownian motion, with the space of
continuous functions on $[0,T]$, the $\sigma$-algebra generated by
the topology of uniform convergence and Wiener measure. The
canonical representation for a pure jump L\'{e}vy process
$(\Omega_{N},\mathcal{F}_{N},\mathbb{P}_{N})$ we use is based on a
fixed partition $(S_{k})_{k\geq 1}$ of $\mathbb{R}-\{0\}$ such that
$0<\nu(S_{k})<\infty, k\geq 1$. Accordingly, it is given by the
product space $\bigotimes_{k\geq
1}(\Omega^{k}_{N},\mathcal{F}^{k}_{N},\mathbb{P}^{k}_{N})$ of spaces
of compound Poisson processes on $[0,T]$ with intensities
$\nu(S_{k})$ and jump size distributions supported on $S_{k}, k\ge
1.$ Since trajectories of compound Poisson processes can be
described by finite families $((t_{1},z_{1}),...,(t_{n},z_{n}))$,
where $(t_{1},...,t_{n})$ denotes the jump times and
$(z_{1},...,z_{n})$ the corresponding sizes of jumps, one can take
$\Omega^{k}_{N}=\bigcup_{n\geq 0}([0,T]\times
(\mathbb{R}-\{0\}))^{n}$, with $([0,T]\times(\mathbb{R}-\{0\})^{0}$
representing an empty sequence, the $\sigma$-algebra
$\mathcal{F}^{k}_{N}=\bigvee_{n\geq
0}\mathcal{B}(([0,T]\times(\mathbb{R}-\{0\}))^{n})$, and the measure
$\mathbb{P}^{k}_{N}$ defined in such a way that for $B=\cup_{n\geq
0}B_{n}, B_{n}\in\mathcal{B}(([0,T]\times(\mathbb{R}-\{0\}))^{n})$,
we have
$$\mathbb{P}^{k}_{N}(B)=e^{-\nu(S_{k})T}\sum_{n=0}^{\infty}\frac{(\nu(S_{k}))^{n}(dt\otimes
\frac{\nu \mathbf{1}\{S_{k}\}}{\nu(S_{k})} )^{\otimes
n}(B_{n})}{n!}.$$ \indent Now consider the finite measure $q$
defined on $[0,T]\times\mathbb{R}$ by
\begin{eqnarray*}
q(E)=\int_{E(0)}dt+\int_{E'}z^{2}\nu(dz)dt,\quad
E\in\mathcal{B}([0,T]\times \mathbb{R}),
\end{eqnarray*}
where $E(0)=\{t\in[0,T]; (t,0)\in E\}$ and $E'=E-E(0)$, and the
random measure $Q$ on $[0,T]\times\mathbb{R}$
\begin{eqnarray*}
Q(E)=\int_{E(0)}dW(t)+\int_{E'}z\tilde{N}(dt,dz),\quad
E\in\mathcal{B}([0,T]\times \mathbb{R}).
\end{eqnarray*}
For $n\in\mathbb{N}$ and a simple function
$h_{n}=\mathbf{1}_{E_{1}\times...\times E_{n}}$, with pairwise
disjoints sets $E_{1},...,E_{n}\in\mathcal{B}([0,T]\times
\mathbb{R})$, a multiple two-parameter integral with respect to the
random measure $Q$
\begin{eqnarray*}
I_{n}(h_{n})=\int_{([0,T]\times\mathbb{R})^{n}}
h((t_{1},z_{1}),...(t_{n},z_{n}))Q(dt_{1},dz_{1})\cdot...\cdot
Q(dt_{n},dz_{n})
\end{eqnarray*}
can be defined as
\begin{eqnarray*}
I_{n}(h_{n})=Q(E_{1})...Q(E_{n}).
\end{eqnarray*}
The integral can be extended to the space
$L^{2}_{T,q,n}(\mathbb{R})$ of product measurable, deterministic
functions $h:([0,T]\times\mathbb{R})^{n}\rightarrow\mathbb{R}$
satisfying
\begin{eqnarray*}
\|h\|^{2}_{L_{T,q,n}^{2}}=\int_{([0,T]\times\mathbb{R})^{n}}
|h_{n}((t_{1},z_{1}),...,(t_{n},z_{n}))|^{2}q(dt_{1},dz_{1})\cdot
...\cdot q(dt_{n},dz_{n})<\infty.
\end{eqnarray*}
The chaotic decomposition property yields that any
$\mathcal{F}$-measurable square integrable random variable $H$ on
the canonical space has a unique representation
\begin{eqnarray}\label{decomposition}
H=\sum_{n=0}^{\infty}I_{n}(h_{n}),\quad \mathbb{P}-a.s.,
\end{eqnarray}
with functions $h_{n}\in L^{2}_{T,q,n}(\mathbb{R})$ that are
symmetric in the $n$ pairs $(t_i,z_i)$, $1\le i\le n$. Moreover,
\begin{eqnarray}\label{square}
\mathbb{E}\big[H^{2}\big]=\sum_{n=0}^{\infty}n!\|h_{n}\|^{2}_{L^{2}_{T,q,n}}.
\end{eqnarray}
In this setting it is possible to study two-parameter annihilation
operators (Malliavin derivatives) and creation operators (Skorokhod
integrals).

\begin{df}\label{def:lspaces}
\begin{itemize}
\item[1.] Let $\mathbb{D}^{1,2}(\mathbb{R})$ denote the space of
$\mathbb{F}$-measurable random variables
$H\in\mathbb{L}^{2}(\mathbb{R})$ with the representation
$H=\sum_{n=0}^{\infty}I_{n}(h_{n})$ satisfying
\begin{eqnarray*}
\sum_{n=1}^{\infty}nn!\|h_{n}\|^{2}_{L_{T,q,n}^{2}}<\infty.
\end{eqnarray*}
\item [2.] The Malliavin derivative $DH:
\Omega\times[0,T]\times\mathbb{R}\rightarrow\mathbb{R}$ of a random
variable $H\in\mathbb{D}^{1,2}(\mathbb{R})$ is a stochastic process
defined by
\begin{eqnarray*}
D_{t,z}H=\sum_{n=1}^{\infty}nI_{n-1}(h_{n}((t,z),\cdot),\,\,\mbox{valid
for}\,\, q-a.e. (t,z)\in[0,T]\times\mathbb{R}, \mathbb{P}-a.s..
\end{eqnarray*}
\item[3.] Let $\mathbb{L}^{1,2}(\mathbb{R})$ denote the space of product measurable and $\mathbb{F}$-adapted processes
$G:\Omega\times[0,T]\times\mathbb{R}\rightarrow\mathbb{R}$
satisfying
\begin{eqnarray*}
&\mathbb{E}\big[\int_{[0,T]\times\mathbb{R}}|G(s,y)|^{2}q(ds,dy)\big]<\infty,\nonumber\\
&G(s,y)\in\mathbb{D}^{1,2}(\mathbb{R}),\,\,\mbox{for}\,\,  q-a.e.  (s,y)\in[0,T]\times\mathbb{R},\nonumber\\
&\mathbb{E}\big[\int_{([0,T]\times\mathbb{R})^{2}}|D_{t,z}G(s,y)|^{2}q(ds,dy)q(dt,dz)\big]<\infty.
\end{eqnarray*}
In terms of the components of the representation of
$G(s,y)=\sum_{n=0}^{\infty}I_{n}(g_{n}((s,y),.)$, for $q$-a.e.
$(s,y)\in[0,T]\times\mathbb{R}$, the above conditions are equivalent
to
$$\sum_{n=1}^{\infty}(n+1)(n+1)!\|\hat{g}_{n}\|^{2}_{L_{T,q,n+1}^{2}}<\infty,$$
where $\hat{g}_n$ denotes the symmetrization of $g_n$ with respect
to all $n+1$
pairs of variables.\\
\noindent The space $\mathbb{L}^{1,2}(\mathbb{R})$ is a Hilbert
space endowed with the norm
\begin{eqnarray*}
\lefteqn{||G||^{2}_{\mathbb{L}^{1,2}}=\mathbb{E}\big[\int_{[0,T]\times\mathbb{R}}|G(s,y)|^{2}q(ds,dy)\big]}\nonumber\\
&&+\mathbb{E}\big[\int_{([0,T]\times\mathbb{R})^{2}}|D_{t,z}G(s,y)|^{2}q(ds,dy)q(dt,dz)\big].
\end{eqnarray*}
\item[4.] The Skorokhod integral with respect to the random measure $Q$ of a process
$G:\Omega\times[0,T]\times\mathbb{R}\rightarrow\mathbb{R}$ with the
representation $G(s,y)=\sum_{n=0}^{\infty}I_{n}(g_{n}((s,y),.)$, for
$q-a.e. (s,y)\in[0,T]\times\mathbb{R}$, satisfying
$$\sum_{n=0}^{\infty}(n+1)!\|\hat{g}_{n}\|^{2}_{L_{T,q,n+1}^{2}}<\infty,$$
is defined as
$$\int_{[0,T]\times\mathbb{R}}G(s,y)Q(ds,dy)=\sum_{n=0}^{\infty}I_{n+1}(\hat{g}_{n}),\quad \mathbb{P}-a.s.$$
\end{itemize}
\end{df}
\indent The following practical rules of differentiation hold.
Consider a random variable $H$ defined on
$\Omega_{W}\times\Omega_{N}$. The derivative $D_{t,0}H$ is with
respect to the Brownian motion component of the L\'{e}vy process,
and we can apply classical Malliavin's calculus for Hilbert
space-valued random variables. If for $\mathbb{P}^{N}$-a.e.
$\omega_{N}\in \Omega_{N}$ the random variable $H(.,\omega_{N})$ is
differentiable in the sense of classical Malliavin's calculus, then
we have the relation
\begin{eqnarray}\label{derbrownian}
D_{t,0}H(\omega_{W},\omega_{N})=D_{t}H(.,\omega_{N})(\omega_{W}),\quad
\lambda-a.e.  t\in[0,T], \mathbb{P}^{W}\times\mathbb{P}^{N}-a.s.,
\end{eqnarray}
where $D_{t}$ denotes the classical Malliavin derivative on the
canonical Brownian space, see Proposition 3.5 in \cite{Sole}. In
order to define $D_{t,z}F$ for $z\neq 0$, which is a derivative with
respect to the pure jump part of the L\'{e}vy process, consider the
following increment quotient operator
\begin{eqnarray}\label{picard}
\Psi_{t,z}H(\omega_{W},\omega_{N})=\frac{H(\omega_{W},\omega^{t,z}_{N})-H(\omega_{W},\omega_{N})}{z},
\end{eqnarray}
where $\omega^{t,z}_{N}$ transforms a family
$\omega_{N}=((t_{1},z_{1}),(t_{2},z_{2}),...))\in\Omega_{N}$ into a
new family
$\omega^{t,z}_{N}=((t,z),(t_{1},z_{1}),(t_{2},z_{2}),...))\in\Omega_{N},$
by adding a jump of size $z$ at time $t$ into the trajectory.
According to Propositions 5.4 and 5.5 in \cite{Sole}, for
$H\in\mathbb{L}^{2}(\mathbb{R})$ such that
$\mathbb{E}\big[\int_{0}^{T}\int_{\mathbb{R}-\{0\}}\big|\Psi_{t,z}H\big|^{2}m(dz)dt
\big]<\infty$ we have the relation
\begin{eqnarray}\label{derpicard}
D_{t,z}H=\Psi_{t,z}H,\,\,\mbox{for}\,\, \lambda\otimes m -a.e.
(t,z)\in[0,T]\times(\mathbb{R}-\{0\}),\mathbb{P}-a.s..
\end{eqnarray}
\indent The operator \eqref{picard} is closely related to Picard's
difference operator, introduced in \cite{Picard}, which is just the
numerator of \eqref{picard}. It is possible to define Malliavin's
derivative for pure jump processes in such a way that it coincides
with Picard's difference operator, see \cite{DIN}. We point out once
again that we adopt the approach of \cite{Sole}, and define multiple
two-parameter integrals with respect to the random measure
$\tilde{M}$ and not with respect to $\tilde{N}$, to obtain
differentiation rules \eqref{derbrownian} and \eqref{derpicard}.
\\
\indent We now discuss some technical problems arising in the next
section in the context of the main theorem of this paper. The
subsequent lemmas are extensions of classical Malliavin
differentiation rules to the setting of the canonical L\'{e}vy
space.

\begin{lem}\label{derexp}
Assume that $H\in\mathbb{D}^{1,2}(\mathbb{R})$. Then, for $0\leq
s\leq T$,
$\mathbb{E}\big[H|\mathcal{F}_{s}\big]\in\mathbb{D}^{1,2}(\mathbb{R})$
and
\begin{eqnarray*}
D_{t,z}\mathbb{E}\big[H|\mathcal{F}_{s}\big]=\mathbb{E}\big[D_{t,z}H|\mathcal{F}_{s}\big]\mathbf{1}\{t\leq
s\},\,\,\mbox{for}\,\, q-a.e. (t,z)\in[0,T]\times\mathbb{R},
\mathbb{P}-a.s..
\end{eqnarray*}
\end{lem}
\Proof The proof is a straightforward extension of the proof of
Proposition 1.2.8 from \cite{N}. Details are left to the
reader.\cbdu \indent We next provide a proof of the commutation of
Lebesgue's integration and variational differentiability, which is
commonly used.
\begin{lem}\label{derint}
Let $G:\Omega\times[0,T]\times\mathbb{R}\to\mathbb{R}$ be a product
measurable and $\mathbb{F}$-adapted process, $\eta$ on
$[0,T]\times\mathbb{R}$ a finite measure, so that the conditions
\begin{eqnarray}\label{monclass}
&\mathbb{E}\big[\int_{[0,T]\times\mathbb{R}}|G(s,y)|^{2}\eta(ds,dy)\big]<\infty,\nonumber\\
&G(s,y)\in\mathbb{D}^{1,2}(\mathbb{R}),\,\,\mbox{for}\,\, \eta-a.e.
(s,y)\in[0,T]\times\mathbb{R},\\
&\mathbb{E}\big[\int_{([0,T]\times\mathbb{R})^{2}}|D_{t,z}G(s,y)|^{2}\eta(ds,dy)q(dt,dz)\big]<\infty.\nonumber
\end{eqnarray}
are satisfied. Then
$\int_{[0,T]\times\mathbb{R}}G(s,y)\eta(ds,dy)\in\mathbb{D}^{1,2}(\mathbb{R})$
and the differentiation rule
\begin{eqnarray*}
D_{t,z}\int_{[0,T]\times\mathbb{R}}G(s,y)\eta(ds,dy)=
\int_{[0,T]\times\mathbb{R}}D_{t,z}G(s,y)\eta(ds,dy)
\end{eqnarray*}
holds for $q$-a.e. $(t,z)\in[0,T]\times\mathbb{R}, \mathbb{P}$-a.s..
\end{lem}
\Proof As for $\eta$-a.e. $(s,y)\in[0,T]\times\mathbb{R}$ the random
variable $G(s,y)$ is $\mathcal{F}_{s}$-measurable and square
integrable, the chaotic decomposition property on the canonical
space \eqref{decomposition} provides a unique representation
\begin{eqnarray*}
G(s,y)=\sum_{n=0}^{\infty}I_{n}(g_{n}((s,y),.)),\quad \eta-a.e.
(s,y)\in[0,T]\times\mathbb{R},\quad \mathbb{P}-a.s..
\end{eqnarray*}
By part 3 of definition \ref{def:lspaces}, the assumptions
\eqref{monclass} yield
\begin{eqnarray}\label{monass}
\int_{[0,T]\times\mathbb{R}}\sum_{n=1}^{\infty}nn!\|g_{n}((s,y),.)\|^{2}_{L^{2}_{T,q,n}}\eta(ds,dy)<\infty.
\end{eqnarray}
For $N\in\mathbb{N}$ let $G^N$ be a measurable version of the
partial sum of the first $N+1$ components given by
\begin{eqnarray*}
G^{N}(s,y)=\sum_{n=0}^{N}I_{n}(g_{n}((s,y),.)),\quad \eta-a.e.
(s,y)\in[0,T]\times\mathbb{R},\quad \mathbb{P}-a.s..
\end{eqnarray*}
We first prove that
$\int_{[0,T]\times\mathbb{R}}G^{N}(s,y)\eta(ds,dy)\in\mathbb{D}^{1,2}(\mathbb{R})$
and the claimed differentiation rule holds. \\
\noindent By canonical extension arguments, a Fubini type property
holds for any single chaos component, and therefore
\begin{eqnarray}\label{ruleN}
\lefteqn{\int_{[0,T]\times\mathbb{R}}G^{N}(s,y)\eta(ds,dy)}\nonumber\\
&=&\int_{[0,T]\times\mathbb{R}}\sum_{n=0}^{N}\int_{([0,T]\times\mathbb{R})^{n}}g_{n}((s,y),(t_{1},z_{1}),...,(t_{n},z_{n}))\nonumber\\
&&\quad \ \quad \cdot Q(dt_{1},dz_{1})...Q(dt_{n},dz_{n})\eta(ds,dy)\nonumber\\
&=&\sum_{n=0}^{N}\int_{([0,T]\times\mathbb{R})^{n}}\int_{[0,T]\times\mathbb{R}}g_{n}((s,y),(t_{1},z_{1}),...,(t_{n},z_{n}))\nonumber\\
&&\quad \ \quad
\cdot\eta(ds,dy)Q(dt_{1},dz_{1})...Q(dt_{n},dz_{n})=\sum_{n=0}^{N}I_{n}(h_{n}):=H^{N},
\end{eqnarray}
with
\begin{eqnarray*}
h_{n}((t_{1},z_{1}),...,(t_{n},z_{n}))=\int_{[0,T]\times\mathbb{R}}g_{n}((s,y),(t_{1},z_{1}),...,(t_{n},z_{n}))\eta(ds,dy)
\end{eqnarray*}
for $(t_{1},z_{1}),...,(t_{n},z_{n})\in ([0,T]\times\mathbb{R})^n.$
Notice that by Cauchy-Schwarz' inequality, finiteness of $\eta$, the
assumption \eqref{monass} and Fubini's theorem we obtain
\begin{eqnarray}\label{fubini}
\sum_{n=1}^{\infty}nn!\|h_{n}\|^{2}_{L^{2}_{T,,q,n}}<\infty.
\end{eqnarray}
For any $N\in\mathbb{N}$, we have that
$H^{N}\in\mathbb{D}^{1,2}(\mathbb{R})$, hence
$\int_{[0,T]\times\mathbb{R}}G^{N}(s,y)\eta(ds,dy)\in\mathbb{D}^{1,2}(\mathbb{R})$,
and, by linearity and definition
\begin{eqnarray*}
\lefteqn{D_{t,z}H^{N}=D_{t,z}\int_{[0,T]\times\mathbb{R}}G^{N}(s,y)\eta(ds,dy)}\\
&=&\sum_{n=1}^{N}n\int_{[0,T]\times\mathbb{R}}g_{n}((s,y),(t,z),(t_{2},z_{2},)...,(t_{n},z_{n}))\eta(ds,dy)Q(dt_{2},dz_{2})...Q(dt_{n},dz_{n})\\
&=&\int_{[0,T]\times\mathbb{R}}\sum_{n=1}^{N}ng_{n}((s,y),(t,z),(t_{2},z_{2},)...,(t_{n},z_{n}))Q(dt_{2},dz_{2})...Q(dt_{n},dz_{n})\eta(ds,dy)\\
&=&\int_{[0,T]\times\mathbb{R}}D_{t,z}G^{N}(s,y)\eta(ds,dy),\quad
q-a.e. (t,z)\in[0,T]\times\mathbb{R}.
\end{eqnarray*}
The differentiation rule is proved for $G^{N}$.\\
\noindent Finally, by \eqref{fubini} we have
\begin{eqnarray*}
\lefteqn{\mathbb{E}\big[|H^{N}-H^{M}|^{2}\big]+\int_{[0,T]\times\mathbb{R}}\mathbb{E}\big[|D_{t,z}H^{N}-D_{t,z}H^{M}|^{2}\big]q(dt,dz)}\\
&&\leq \sum_{n=N+1}^{M}nn!\|h_{n}\|^{2}_{L^{2}_{T,q,n}}\rightarrow
0, \quad N,M\rightarrow\infty.
\end{eqnarray*}
By closeability of the operator $D$ we conclude that the unique
limit $H$ is Malliavin differentiable. The convergences
$G^{N}\rightarrow G$ $\mathbb{P}\otimes\eta$-a.e. and $D
G^{N}\rightarrow D G$ $\mathbb{P}\otimes\eta\otimes q$-a.e. together
with Lebesgue's dominated convergence theorem, justified by the
first and third assumption in \eqref{monclass}, give
\begin{eqnarray*}
\lefteqn{\mathbb{E}\big[|\int_{[0,T]\times\mathbb{R}}G^{N}(s,y)\eta(ds,dy)-\int_{[0,T]\times\mathbb{R}}G(s,y)\eta(ds,dy)|^{2}\big]}\\
&&+\int_{[0,T]\times\mathbb{R}}
\mathbb{E}\big[|\int_{[0,T]\times\mathbb{R}}D_{t,z}G^{N}(s,y)\eta(ds,dy)-\int_{[0,T]\times\mathbb{R}}D_{t,z}G(s,y)\eta(ds,dy)|^{2}\big]q(dt,dz)\rightarrow
0.
\end{eqnarray*}
This implies the claimed equation.\cbdu

We finally discuss the commutation relation of the Skorokhod
stochastic integral with the variational derivative.

\begin{lem}\label{derl}
Assume that $G:\Omega\times[0,T]\times\mathbb{R}\to\mathbb{R}$ is a
predictable process and
$\mathbb{E}\big[\int_{[0,T]\times\mathbb{R}}|G(s,y)|^{2}q(ds,dy)\big]<\infty$
holds. Then
\begin{center}
$G\in\mathbb{L}^{1,2}(\mathbb{R})$ if and only if
$\int_{[0,T]\times\mathbb{R}}G(s,y)Q(ds,dy)\in\mathbb{D}^{1,2}(\mathbb{R})$.
\end{center}
Moreover, if
$\int_{[0,T]\times\mathbb{R}}G(s,y)Q(ds,dy)\in\mathbb{D}^{1,2}(\mathbb{R})$
then, for $q$-a.e $(t,z)\in[0,T]\times\mathbb{R}$,
\begin{eqnarray*}
D_{t,z}\int_{[0,T]\times\mathbb{R}}G(s,y)Q(ds,dy)
=G(t,z)+\int_{[0,T]\times\mathbb{R}}D_{t,z}G(s,y)Q(ds,dy),\quad
\mathbb{P}-a.s.,
\end{eqnarray*}
and $\int_{[0,T]\times\mathbb{R}}D_{t,z}G(s,y)Q(ds,dy)$ is a
stochastic integral in It\^{o} sense.
\end{lem}
\Proof By square integrability of $G$, for $q$-a.e
$(s,y)\in[0,T]\times\mathbb{R}$, the chaotic decomposition property
yields the unique representation
$G(s,y)=\sum_{n=0}^{\infty}I_{n}(g_{n}((s,y),.)$, $g_{n}\in
L^{2}_{T,q,n+1},n\geq 0$. Square integrability and predictability of
$G$ implies that the stochastic integral
$\int_{[0,T]\times\mathbb{R}}G(s,y)Q(ds,dy)$ is well-defined in the
It\^{o} sense and the Skorokhod integral, which coincides under the
given assumptions with the It\^{o} integral (see Theorem 6.1 in
\cite{Sole}) can be defined by the series expansion
$\int_{[0,T]\times\mathbb{R}}G(s,y)Q(ds,dy)=\sum_{n=0}^{\infty}I_{n+1}(\hat{g}_{n})$
according to Definition 3.1.4. The Skorokhod integral is Malliavin
differentiable if and only if
$\sum_{n=1}^{\infty}(n+1)(n+1)!\|\hat{g}_{n}\|^{2}_{L^{2}_{T,q,n+1}}<\infty$,
see Definition 3.1.2. This series converges if and only if
$G\in\mathbb{L}^{1,2}(\mathbb{R})$, by Definition 3.1.3. \\
\noindent Following Section 6 in \cite{Sole}, we can conclude that
the required differentiation rule holds. To prove that the integral
$\int_{[0,T]\times\mathbb{R}}D_{t,z}G(s,y)Q(ds,dy)$ is well-defined
in the It\^{o} sense, it is sufficient to show that the integrand
$(\omega,s,y)\mapsto D_{t,z}G(s,y)(\omega)$ is a predictable mapping
on $\Omega\times[0,T]\times\mathbb{R}$, as square integrability is
already satisfied by $G\in\mathbb{L}^{1,2}(\mathbb{R})$. For
$q$-a.e. $(s,y)\in[0,T]\times\mathbb{R}$, predictability of $G$
implies that
\begin{eqnarray*}
G(s,y)=\sum_{n=0}^{\infty}I_{n}(g_{n}((s,y),.)=\sum_{n=0}^{\infty}I_{n}(g_{n}((s,y),.)\mathbf{1}^{\otimes
n}_{[0,s)}(.)),\quad \mathbb{P}-a.s.,
\end{eqnarray*}
and applying Definition 3.1.2 of the Malliavin derivative yields
\begin{eqnarray*}
\lefteqn{D_{t,z}G(s,y)=\sum_{n=0}^{\infty}nI_{n-1}(g_{n}((s,y),(t,z),.)\mathbf{1}^{\otimes
n}_{[0,s)}((t,z),.)),}\nonumber\\
&&\mbox{for}\,\,q\otimes q-a.e. ((t,z),
(s,y))\in([0,T]\times\mathbb{R})^{2}, \mathbb{P}-a.s.,
\end{eqnarray*}
from which the required predictability of the integrand follows. As
a by-product, let us note that $(\omega,s,y,t,z)\mapsto
D_{t,z}G(s,y)(\omega)$ is jointly measurable. \cbdu

\section{Variational differentiability of a solution}
\indent The main goal of this paper is to investigate Malliavin's
differentiability of a solution of a backward stochastic
differential equation with a time delayed generator. In this
section, additionally to \textbf{(A1)}-\textbf{(A5)}, we assume that
\begin{description}
    \item[(A6)] the generator $f$ is of the following form
    \begin{eqnarray*}
    \lefteqn{f(t,y_{t},z_{t},u_{t})}\\
    &:=&f\Big(\omega,t,\int_{-T}^{0}y(t+v)\alpha(dv),\\
    &&\quad \ \quad \int_{-T}^{0}z(t+v)\alpha(dv)\int_{-T}^{0}\int_{\mathbb{R}-\{0\}}u(t+v,z)m(dz)\alpha(dv)\Big),
    \end{eqnarray*}
    with a product measurable function
    $f:\Omega\times[0,T]\times\mathbb{R}\times\mathbb{R}\times\mathbb{R}\rightarrow\mathbb{R}$,
    which is Lipschitz continuous in the last three variables for
    $\mathbb{P}\otimes\lambda$-a.e.
    $(\omega,t)\in\Omega\times[0,T]$, more precisely the generator satisfies \textbf{(A3)} with
    the same constant $K$,
  \item[(A7)] the terminal value is Malliavin differentiable, i.e.
    $\xi\in\mathbb{D}^{1,2}(\mathbb{R})$, and
    \begin{eqnarray*}
    \mathbb{E}\big[\int_{[0,T]\times\mathbb{R}}|D_{s,z}\xi|^{2}q(ds,dz)\big]<\infty,\nonumber\\
    \lim_{\epsilon\downarrow
    0}\mathbb{E}\big[\int_{0}^{T}\int_{|z|\leq\epsilon}|D_{s,z}\xi|^{2}m(dz)ds\big]=0,
    \end{eqnarray*}
  \item[(A8)] for $\mathbb{P}\otimes\lambda$-a.e. $(\omega,t)\in\Omega\times[0,T]$, the mapping $(y,z,u)\mapsto f(\omega,t,y,z,u)$
  is continuously differentiable in $(y,z,u)$, with uniformly bounded and continuous partial derivatives $f_{y}, f_{z},
    f_{u}$; we assume $f_{y}(\omega,t,.,.,.)=f_{z}(\omega,t,.,.,.)=f_{u}(\omega,t,.,.,.)=0$ for
    $\omega\in\Omega,t<0$;
    \item[(A9)] for
    $(t,y,z,u)\in[0,T]\times\mathbb{R}\times\mathbb{R}\times\mathbb{R}$
    we have
    $f(\cdot,t,y,z,u)\in\mathbb{D}^{1,2}(\mathbb{R})$ and
    \begin{eqnarray*}
    &&\mathbb{E}\big[\int_{[0,T]\times\mathbb{R}}\int_{0}^{T}\big|D_{s,z}f(\cdot,t,0,0,0)\big|^{2}dt\,q(ds,dy)\big]<\infty,\\
    \lefteqn{\big|D_{s,z}f(\omega,t,\hat{y},\hat{z},\hat{u})-\big|D_{s,z}f(\omega,t,\tilde{y},\tilde{z},\tilde{u}\big)\big|}\\
    &&\leq L\big(|\hat{y}-\tilde{y}|+|\hat{z}-\tilde{z}|+|\hat{u}-\tilde{u}|\big),\quad
     \end{eqnarray*}
    $(s,z)\in[0,T]\times\mathbb{R}$, $(\hat{y},\hat{z},\hat{u})\in\mathbb{R}\times\mathbb{R}\times\mathbb{R},
    (\tilde{y},\tilde{z},\tilde{u})\in\mathbb{R}\times\mathbb{R}\times\mathbb{R}$, for $\mathbb{P}\otimes\lambda$-a.e.
    $(\omega,t)\in\Omega\times[0,T].$
\end{description}
The assumptions \textbf{(A7)}-\textbf{(A9)} are classical when
dealing with Malliavin's differentiability, see Proposition 5.3 in
\cite{K} or Theorem 3.3.1 in \cite{I1} in the case of BSDEs driven
by Brownian motions. We also remark that the generator in
\textbf{(A6)} depends on $\int_{-T}^{0}\int_{\mathbb{R-}\{0\}}u(t+v,
z)m(dz)\alpha(dv)$, which corresponds to a standard form of
dependence appearing in BSDE without delays and with jumps,
see Proposition 2.6 and Remark 2.7 in \cite{BBP}.\\
\indent We can state our main theorem.
\begin{thm}\label{thmmalliavin}
Assume that \textbf{(A1)}-\textbf{(A9)} hold and that time horizon
$T$ and Lipschitz constant $K$ of the generator $f$ are sufficiently
small, such that for some $\beta>0$
\begin{eqnarray*}
\delta:=\delta(T,K,\beta,\alpha)=(8T+\frac{1}{\beta})K\int_{-T}^{0}e^{-\beta
v}\alpha(dv)\max\{1,T\}<1.
\end{eqnarray*}
\noindent 1. There exists a unique solution
$(Y,Z,U)\in\mathbb{S}^{2}(\mathbb{R})\times
\mathbb{H}^{2}(\mathbb{R})\times\mathbb{H}_{m}^{2}(\mathbb{R})$ of
the time delayed BSDE
\begin{eqnarray}\label{bsdem}
\lefteqn{Y(t)=\xi}\nonumber\\
&&+\int_{t}^{T}f\Big(\omega,r,\int_{-T}^{0}Y(r+v)\alpha(dv),\nonumber\\
&&\quad \ \quad \int_{-T}^{0}Z(r+v)\alpha(dv),\int_{-T}^{0}\int_{\mathbb{R}-\{0\}}U(r+v,z)m(dz)\alpha(dv)\Big)dr\nonumber\\
&&-\int_{t}^{T}Z(r)dW(r)-\int_{t}^{T}\int_{\mathbb{R}-\{0\}}U(r,y)\tilde{M}(dr,dy),\quad
0\leq t\leq T.
\end{eqnarray}
\noindent 2. There exists a unique solution
$(Y^{s,0},Z^{s,0},U^{s,0})\in\mathbb{S}^{2}(\mathbb{R})\times
\mathbb{H}^{2}(\mathbb{R})\times\mathbb{H}_{m}^{2}(\mathbb{R})$ of
the time delayed BSDE
\begin{eqnarray}\label{eq0}
Y^{s,0}(t)&=&D_{s,0}\xi+\int_{t}^{T}f^{s,0}(r)dr-\int_{t}^{T}Z^{s,0}(r)dW(r)\nonumber\\
&&-
\int_{t}^{T}\int_{\mathbb{R}-\{0\}}U^{s,0}(r,y)\tilde{M}(dr,dy),\quad
0\leq s\leq t\leq T,
\end{eqnarray}
with the generator
\begin{eqnarray}\label{op0}
\lefteqn{f^{s,0}(r)}\nonumber\\
&=&D_{t,0}f\Big(\omega,r,\int_{-T}^{0}Y(r+v)\alpha(dv),\nonumber\\
&&\quad \quad \int_{-T}^{0}Z(r+v)\alpha(dv),\int_{-T}^{0}\int_{\mathbb{R}-\{0\}}U(r+v,y)m(dy)\alpha(dv)\Big)\nonumber\\
&+&f_{y}\Big(\omega,r,\int_{-T}^{0}Y(r+v)\alpha(dv),\nonumber\\
&&\quad \ \quad \int_{-T}^{0}Z(r+v)\alpha(dv),\int_{-T}^{0}\int_{\mathbb{R}-\{0\}}U(r+v,y)m(dy)\alpha(dv)\Big)\nonumber\\
&&\cdot\int_{-T}^{0}Y^{s,0}(r+v)\alpha(dv)\nonumber\\
&+&f_{z}\Big(\omega,r,\int_{-T}^{0}Y(r+v)\alpha(dv),\nonumber\\
&&\quad \ \quad \int_{-T}^{0}Z(r+v)\alpha(dv),\int_{-T}^{0}\int_{\mathbb{R}-\{0\}}U(r+v,y)m(dy)\alpha(dv)\Big)\nonumber\\
&&\cdot\int_{-T}^{0}Z^{s,0}(r+v)\alpha(dv)\nonumber\\
&+&f_{u}\Big(\omega,r,\int_{-T}^{0}Y(r+v)\alpha(dv),\nonumber\\
&&\quad \quad \int_{-T}^{0}Z(r+v)\alpha(dv),\int_{-T}^{0}\int_{\mathbb{R}-\{0\}}U(r+v,y)m(dy)\alpha(dv)\Big)\nonumber\\
&&\cdot\int_{-T}^{0}\int_{\mathbb{R}-\{0\}}U^{s,0}(r+v,y)m(dy)\alpha(dv).
\end{eqnarray}
\noindent 3. There exists a unique solution
$(Y^{s,z},Z^{s,z},U^{s,z})\in\mathbb{S}^{2}(\mathbb{R})\times
\mathbb{H}^{2}(\mathbb{R})\times\mathbb{H}_{m}^{2}(\mathbb{R})$ of
the time delayed BSDE
\begin{eqnarray}\label{eqz}
Y^{s,z}(t)&=&D_{s,z}\xi+\int_{t}^{T}f^{s,z}(r)dr-\int_{t}^{T}Z^{s,z}(r)dW(r)\nonumber\\
&&-
\int_{t}^{T}\int_{\mathbb{R}-\{0\}}U^{s,z}(r,y)\tilde{M}(dr,dy),\quad
0\leq s\leq t\leq T,z\neq 0,
\end{eqnarray}
with the generator
\begin{eqnarray}\label{opz}
\lefteqn{f^{s,z}(r)}\nonumber\\
&=&\left\{f\Big(\omega^{s,z},r,z\int_{-T}^{0}Y^{s,z}(r+v)\alpha(dv)+\int_{-T}^{0}Y(r+v)\alpha(dv),\right.\nonumber\\
&&\quad  z\int_{-T}^{0}Z^{s,z}(r+v)\alpha(dv)+\int_{-T}^{0}Z(r+v)\alpha(dv),\nonumber\\
&&\quad  z\int_{-T}^{0}\int_{\mathbb{R}-\{0\}}U^{s,z}(r+v,y)m(dy)\alpha(dv)+\int_{-T}^{0}\int_{\mathbb{R}-\{0\}}U(r+v,x)m(dy)\alpha(dv)\Big)\nonumber\\
&&-f\Big(\omega,r,\int_{-T}^{0}Y(r+v)\alpha(dv),\nonumber\\
&&\left.\quad \  \int_{-T}^{0}Z(r+v)\alpha(dv),
\int_{-T}^{0}\int_{\mathbb{R}-\{0\}}U(r+v,y)m(dy)\alpha(dv)\Big)\right\}/z,
\end{eqnarray}
where we set
\begin{eqnarray}\label{zero}
Y^{s,z}(t)=Z^{s,z}(t)=U^{s,z}(t,y)=0, \quad
(y,z)\in(\mathbb{R}-\{0\})\times\mathbb{R}, \mathbb{P}-a.s., t<s\leq
T.
\end{eqnarray}
Then
$(Y,Z,U)\in\mathbb{L}^{1,2}(\mathbb{R})\times\mathbb{L}^{1,2}(\mathbb{R})\times\mathbb{L}^{1,2}(\mathbb{R})$
and\\ $(Y^{s,z}(t),Z^{s,z}(t),U^{s,z}(t,y))_{0\leq s,t\leq
T,(y,z)\in(\mathbb{R}-\{0\})\times\mathbb{R}}$ is a version of\\
$(D_{s,z}Y(t),D_{s,z}Z(t),D_{s,z}U(t,y))_{0\leq s,t\leq
T,(y,z)\in(\mathbb{R}-\{0\})\times\mathbb{R}}$.
\end{thm}

We recall that $D_{t,0}f(r,.,.,.,.)$ appearing as the first term in
\eqref{op0} is the Malliavin derivative of $f$ with respect to
$\omega$, whereas $\omega^{s,z}$ appearing in \eqref{opz} is defined
in \eqref{picard}.

\Proof We follow the idea of the proofs of Proposition 5.3 in
\cite{K}, or Theorem 3.3.1 in \cite{I1}. Let us denote by $C$
a finite constant which may change from line to line.\\
\noindent Step 1) Given $\beta>0$, we prove existence of unique
solutions of the equations \eqref{bsdem}, \eqref{eq0} and
\eqref{eqz}
for a time horizon $T$ and a Lipschitz constant $K$ fulfilling $\delta(T,K,\beta,\alpha)<1$.\\
\noindent The existence of a unique solution
$(Y,Z,U)\in\mathbb{S}^{2}(\mathbb{R})\times
\mathbb{H}^{2}(\mathbb{R})\times\mathbb{H}_{m}^{2}(\mathbb{R})$ of
\eqref{bsdem} follows from Theorem \ref{thmain}, since the
assumptions \textbf{(A1)}-\textbf{(A5)} are satisfied. Under the
additional assumptions \textbf{(A7)}-\textbf{(A9)}, the time delayed
BSDE \eqref{eq0} and \eqref{eqz}, with the generators \eqref{op0}
resp. \eqref{opz}, fulfill the conditions of Theorem \ref{thmain}.
In particular the corresponding generators are Lipschitz continuous
with the same Lipschitz constant $K$ that the generator $f$
possesses. It is easy to see that the generators \eqref{op0} and
\eqref{opz} have the same Lipschitz constant $K$ in the sense of
\textbf{(A3)}. Hence, $\delta(T,K,\beta,\alpha)<1$ holds
simultaneously for all BSDEs \eqref{bsdem}, \eqref{eq0} and
\eqref{eqz} and we  conclude that for $q$-a.e.
$(s,z)\in[0,T]\times\mathbb{R}$ there exists a unique solution
$(Y^{s,z},Z^{s,z},U^{s,z})\in\mathbb{S}^{2}(\mathbb{R})\times
\mathbb{H}^{2}(\mathbb{R})\times\mathbb{H}_{m}^{2}(\mathbb{R})$ of
\eqref{eq0} or \eqref{eqz} satisfying \eqref{zero}.\\

\noindent Step 2) \noindent Consider a sequence
$(Y^{n},Z^{n},U^{n})_{n\in\mathbb{N}}$, constructed by Picard iteration scheme, which converges to $(Y,Z,U)$. In this step we show that
$(Y^{n},Z^{n},U^{n})\in\mathbb{L}^{1,2}(\mathbb{R})\times\mathbb{L}^{1,2}(\mathbb{R})\times\mathbb{L}^{1,2}(\mathbb{R})$
implies
$(Y^{n+1},Z^{n+1},U^{n+1})\in\mathbb{L}^{1,2}(\mathbb{R})\times\mathbb{L}^{1,2}(\mathbb{R})
\times\mathbb{L}^{1,2}(\mathbb{R})$, and that from
$\mathbb{E}\big[\int_{[0,T]}\sup_{t\in[0,T]}|D_{s,z}Y^{n}(t)|^{2}q(ds,dz)\big]<\infty$
we can as well deduce that
$\mathbb{E}\big[\int_{[0,T]\times\mathbb{R}}\sup_{t\in[0,T]}|D_{s,z}Y^{n+1}(t)|^{2}q(ds,dz)\big]<\infty$.\\

\noindent For that purpose, we study the iterations
\begin{eqnarray}\label{bsdenn}
\lefteqn{Y^{n+1}(t)=\xi+\int_{t}^{T}f^{n}(r)dr}\nonumber\\
&&-\int_{t}^{T}Z^{n+1}(r)dW(r)-\int_{t}^{T}\int_{\mathbb{R}-\{0\}}U^{n+1}(r,y)\tilde{M}(dr,dy),\quad
0\leq t\leq T,
\end{eqnarray}
where we denote
\begin{eqnarray*}
\lefteqn{f^{n}(r)}\nonumber\\
&=&f(r,\int_{-T}^{0}Y^{n}(r+v)\alpha(dv),\int_{-T}^{0}Z^{n}(r+v)\alpha(dv),
\int_{-T}^{0}\int_{\mathbb{R}-\{0\}}U^{n}(r+v,y)m(dy)\alpha(dv)).
\end{eqnarray*}
We first establish Malliavin's differentiability  of
$\int_{t}^{T}f^{n}(r)dr$ by applying Lemma \ref{derint}. Notice that
$Y^{n}(t)\in \mathbb{D}^{1,2}(\mathbb{R})$, for $\lambda$-a.e.
$t\in[-T,T]$. Similarly to \eqref{estf}, we can derive
\begin{eqnarray*}
&&\int_{0}^{T}\mathbb{E}\big[\int_{-T}^{0}|Y^{n}(r+v)|^{2}\alpha(dv)\big]dr=
\mathbb{E}\big[\int_{-T}^{0}\int_{0}^{T}|Y^{n}(r+v)|^{2}dr\alpha(dv)\big]\nonumber\\
&& =
\mathbb{E}\big[\int_{-T}^{0}\int_{v}^{T+v}|Y^{n}(w)|^{2}dw\alpha(dv)\big]\leq
T\mathbb{E}\big[\sup_{w\in[0,T]}|Y^{n}(w)|^{2})\big]<\infty
\end{eqnarray*}
together with
\begin{eqnarray*}
\lefteqn{\int_{0}^{T}\mathbb{E}\big[\int_{[0,T]\times\mathbb{R}}\int_{-T}^{0}|D_{s,z}Y^{n}(r+v)|^{2}\alpha(dv)q(ds,dz)\big]dr}\nonumber\\
&\leq&T\mathbb{E}\big[\int_{[0,T]\times
\mathbb{R}}\sup_{w\in[0,T]}|D_{s,z}Y^{n}(w)|^{2}q(ds,dz)\big]<\infty.
\end{eqnarray*}
This provides the assumptions of Lemma \ref{derint}, and for
$\lambda$-a.e. $r\in[0,T]$ we have
$\int_{-T}^{0}Y^{n}(r+v)\alpha(dv)\in\mathbb{D}^{1,2}(\mathbb{R})$,
and furthermore
\begin{eqnarray*}
D_{s,z}\int_{-T}^{0}Y^{n}(r+v)\alpha(dv)=\int_{-T}^{0}D_{s,z}Y^{n}(r+v)\alpha(dv),\quad
 \mathbb{P}-a.s.,
\end{eqnarray*}
for $q\otimes\lambda$-a.e.
$(s,z,r)\in[0,T]\times\mathbb{R}\times[0,T]$. In an analogous way we
derive
\begin{eqnarray*}
&&D_{s,z}\int_{-T}^{0}Z^{n}(r+v)\alpha(dv)=\int_{-T}^{0}D_{s,z}Z^{n}(r+v)\alpha(dv),\\
&&D_{s,z}\int_{-T}^{0}\int_{\mathbb{R}-\{0\}}U^{n}(r+v,y)m(dy)\alpha(dv)=\int_{-T}^{0}\int_{\mathbb{R}-\{0\}}D_{s,z}U^{n}(r+v,y)m(dy)\alpha(dv),
\end{eqnarray*}
holds $\mathbb{P}$-a.s. for $q\otimes\lambda$-a.e.
$(s,z,r)\in[0,T]\times\mathbb{R}\times[0,T]$. We claim that for
$\lambda$-a.e. $r\in[0,T]$ the random variable
$f^{n}(r)\in\mathbb{D}^{1,2}(\mathbb{R})$ and for
$q\otimes\lambda$-a.e. $(s,z,r)\in[0,T]\times\mathbb{R}\times[0,T]$
we have
\begin{eqnarray}\label{op0n}
\lefteqn{D_{s,0}f^{n}(r)}\nonumber\\
&=&D_{t,0}f\Big(\cdot,r,\int_{-T}^{0}Y^{n}(r+v)\alpha(dv),\nonumber\\
&&\quad \ \quad\int_{-T}^{0}Z^{n}(r+v)\alpha(dv),
\int_{-T}^{0}\int_{\mathbb{R}-\{0\}}U^{n}(r+v,y)m(dy)\alpha(dv)\Big)\nonumber\\
&+&f_{y}\Big(\cdot,r,\int_{-T}^{0}Y^{n}(r+v)\alpha(dv),\nonumber\\
&&\quad \ \quad \int_{-T}^{0}Z^{n}(r+v)\alpha(dv), \int_{-T}^{0}\int_{\mathbb{R}-\{0\}}U^{n}(r+v,y)m(dy)\alpha(dv)\Big)\nonumber\\
&&\cdot\int_{-T}^{0}D_{s,0}Y^{n}(r+v)\alpha(dv)\nonumber\\
&+&f_{z}\Big(\cdot,r,\int_{-T}^{0}Y^{n}(r+v)\alpha(dv),\nonumber\\
&&\quad \ \quad\int_{-T}^{0}Z^{n}(r+v)\alpha(dv),
\int_{-T}^{0}\int_{\mathbb{R}-\{0\}}U^{n}(r+v,y)m(dy)\alpha(dv)\Big)\nonumber\\
&&\cdot\int_{-T}^{0}D_{s,0}Z^{n}(r+v)\alpha(dv)\nonumber\\
&+&f_{u}\Big(\cdot,r,\int_{-T}^{0}Y^{n}(r+v)\alpha(dv),\nonumber\\
&&\quad \ \quad \int_{-T}^{0}Z^{n}(r+v)\alpha(dv),
\int_{-T}^{0}\int_{\mathbb{R}-\{0\}}U^{n}(r+v,y)m(dy)\alpha(dv)\Big)\nonumber\\
&&\cdot\int_{-T}^{0}\int_{\mathbb{R}-\{0\}}D_{s,0}U^{n}(r+v,y)m(dy)\alpha(dv),
\end{eqnarray}
and for $z\neq 0$
\begin{eqnarray}\label{opzn}
\lefteqn{D_{s,z}f^{n}(r)}\nonumber\\
&=&\left\{f\Big(\cdot^{s,z},r,z\int_{-T}^{0}D_{s,z}Y^{n}(r+v)\alpha(dv)+\int_{-T}^{0}Y^{n}(r+v)\alpha(dv),\right.\nonumber\\
&&\quad  z\int_{-T}^{0}D_{s,z}Z^{n}(r+v)\alpha(dv)+\int_{-T}^{0}Z^{n}(r+v)\alpha(dv),\nonumber\\
&&\quad  z\int_{-T}^{0}\int_{\mathbb{R}-\{0\}}D_{s,z}U^{n}(r+v,y)m(dy)\alpha(dv)+\int_{-T}^{0}\int_{\mathbb{R}-\{0\}}U^{n}(r+v,x)m(dy)\alpha(dv)\Big)\nonumber\\
&&-f\Big(\cdot,r,\int_{-T}^{0}Y^{n}(r+v)\alpha(dv),\nonumber\\
&&\left.\quad \ \int_{-T}^{0}Z^{n}(r+v)\alpha(dv),
\int_{-T}^{0}\int_{\mathbb{R}-\{0\}}U^{n}(r+v,y)m(dy)\alpha(dv)\Big)\right\}/z.
\end{eqnarray}
The derivative \eqref{op0n} follows from the chain rule for the
operator $D_{s,0}$, as for Theorem 2 in \cite{P}, whereas
\eqref{opzn} follows from Proposition 5.5 in \cite{Sole} provided
that
\begin{eqnarray*}
&&\mathbb{E}\big[\int_{0}^{T}|f^{n}(r)|^{2}dr\big]<\infty,\\
&&\mathbb{E}\big[\int_{0}^{T}\int_{0}^{T}\int_{\mathbb{R}-\{0\}}|D_{s,z}f^{n}(r)|^{2}m(dz)dsdr\big]<\infty,
\end{eqnarray*}
hold. The finiteness of the first integral is obvious. The second
integral can be shown to be finite by applying the Lipschitz
continuity of the generator \textbf{(A3)}, the Lipschitz continuity
of the derivative of the function $f$ with respect to $\omega$ and
its square integrability \textbf{(A9)}, as well as the assumption
$(Y^{n},Z^{n},U^{n})\in\mathbb{L}^{1,2}(\mathbb{R})\times\mathbb{L}^{1,2}(\mathbb{R})\times\mathbb{L}^{1,2}(\mathbb{R})$.
Moreover,
\begin{eqnarray*}
\mathbb{E}\big[\int_{0}^{T}\int_{0}^{T}|D_{s,0}f^{n}(r)|^{2}drds\big]<\infty,
\end{eqnarray*}
and by Lemma \ref{derint} again we derive that for $0\leq t\leq T$
we have $\xi+\int_{t}^{T}f^{n}(r)dr\in\mathbb{D}^{1,2}(\mathbb{R})$
with Malliavin derivative
\begin{eqnarray}\label{mf}
D_{s,z}\xi+\int_{t}^{T}D_{s,z}f^{n}(r)dr,\quad q-a.e.
(s,z)\in[0,T]\times \mathbb{R},
\end{eqnarray}
where $D_{s,z}f^{n}$ is defined in \eqref{op0n} and \eqref{opzn}. If
we combine this result with Lemma \ref{derexp}, we can conclude
\begin{eqnarray*}
Y^{n+1}(t)=\mathbb{E}\big[\xi+\int_{t}^{T}f^{n}(r)dr|\mathcal{F}_{t}\big]\in\mathbb{D}^{1,2}(\mathbb{R}),\quad
0\leq t\leq T,
\end{eqnarray*}
and from the equation\eqref{bsdenn} we derive
\begin{eqnarray}\label{skor1}
\int_{t}^{T}Z^{n+1}(r)dW(r)\in\mathbb{D}^{1,2}(\mathbb{R}),\quad
0\leq t\leq T,
\end{eqnarray}
and
\begin{eqnarray}\label{skor2}
\int_{t}^{T}\int_{\mathbb{R}-\{0\}}U^{n+1}(r,y)\tilde{M}(dr,dy)\in\mathbb{D}^{1,2}(\mathbb{R}),\quad
0\leq t\leq T.
\end{eqnarray}
Therefore Lemma \ref{derl} yields
$(Z^{n+1},U^{n+1})\in\mathbb{L}^{1,2}(\mathbb{R})\times\mathbb{L}^{1,2}(\mathbb{R})$.\\
\noindent This allows us to differentiate the recursive equation
\eqref{bsdenn} and obtain for
$q$-a.e.$(s,z)\in[0,T]\times\mathbb{R}$
\begin{eqnarray}\label{eqztn1}
D_{s,z}Y^{n+1}(t)&=&D_{s,z}\xi+\int_{t}^{T}D_{s,z}f^{n}(r)dr-\int_{t}^{T}D_{s,z}Z^{n+1}(r)dW(r)\nonumber\\
&&-\int_{t}^{T}\int_{\mathbb{R}-\{0\}}D_{s,z}U^{n+1}(r,y)\tilde{M}(dr,dy),\quad
s\leq t\leq T,
\end{eqnarray}
and \begin{eqnarray}\label{eqztn2}
D_{s,z}Y^{n+1}(t)=D_{s,z}Z^{n+1}(t)=D_{s,z}U^{n+1}(t,y)=0, \quad
t<s, y\in(\mathbb{R}-\{0\}).
\end{eqnarray}
Note that the time delayed BSDE \eqref{eqztn1} with generator
\eqref{op0n} or \eqref{opzn} fulfills the assumptions of Theorem
\ref{thmain} with zero corresponding Lipschitz constant. We conclude
that for $q$-a.e.$(s,z)\in[0,T]\times\mathbb{R}$ there exists a
unique solution $(D_{s,z}Y^{n+1},D_{s,z}Z^{n+1},D_{s,z}U^{n+1})\in
\mathbb{S}^{2}\times\mathbb{H}^{2}\times\mathbb{H}^{2}_{m}$ of
\eqref{eqztn1} satisfying \eqref{eqztn2}. By applying Lemma
\ref{priori}, with $\tilde{\xi}=0$ and $\tilde{f}=f$, together with
the estimate \eqref{contr2}, with $\delta =
\delta(T,K,\beta,\alpha)<1$ we derive the inequality
\begin{eqnarray}\label{recconv1}
\lefteqn{\big\|D_{s,z}Y^{n+1}\big\|^{2}_{\mathbb{S}^{2}}+\big\|D_{s,z}Z^{n+1}\big\|^{2}_{\mathbb{H}^{2}}
+\big\|D_{s,z}U^{n+1}\big\|^{2}_{\mathbb{H}^{2}_{m}}}\nonumber\\
&\leq& 9e^{\beta
T}\mathbb{E}\big[|D_{s,z}\xi|^{2}\big]+\delta\big(\big\|D_{s,z}Y^{n}\big\|^{2}_{\mathbb{S}^{2}}+\big\|D_{s,z}Z^{n}\big\|^{2}_{\mathbb{H}^{2}}+
\big\|D_{s,z}U^{n}\big\|^{2}_{\mathbb{H}^{2}_{m}}\big).
\end{eqnarray}
This in turn yields
$\mathbb{E}\big[\int_{[0,T]\times\mathbb{R}}\sup_{t\in[0,T]}|D_{s,z}Y^{n+1}(t)|^{2}q(ds,dz)\big]<\infty,$
and in particular, $Y^{n+1}\in\mathbb{L}^{1,2}(\mathbb{R})$.\\

\noindent Step 3) We establish the integrability of the solution
$Y^{s,z}(t),Z^{s,z}(t),U^{s,z}(t,y)$ with respect to the product
measure $q$ on ($[0,T]\times\mathbb{R})^{2}$.\\
\noindent Take $(s,z)\in[0,T]\times \mathbb{R}$. Consider the unique
solution
$(Y^{s,z},Z^{s,z},U^{s,z})\in\mathbb{S}^{2}(\mathbb{R})\times
\mathbb{H}^{2}(\mathbb{R})\times\mathbb{H}_{m}^{2}(\mathbb{R})$ of
the equation \eqref{eq0} or \eqref{eqz}. Lemma \ref{priori}, with
$\tilde{\xi}=0$ and $\tilde{f}=f$ together with the estimates
\eqref{contr2} and \eqref{contr} yield the inequality
\begin{eqnarray*}
\lefteqn{\big\|Y^{s,z}\big\|^{2}_{\mathbb{S}^{2}}+\big\|Z^{s,z}\big\|^{2}_{\mathbb{H}^{2}}+
\big\|U^{s,z}\big\|^{2}_{\mathbb{H}^{2}_{m}}}\nonumber\\
&\leq& 9e^{\beta
T}\mathbb{E}\big[|D_{s,z}\xi|^{2}\big]+\delta\big(\big\|Y^{s,z}\big\|^{2}_{\mathbb{S}^{2}}+\big\|Z^{s,z}\big\|^{2}_{\mathbb{H}^{2}}+
\big\|U^{s,z}\big\|^{2}_{\mathbb{H}^{2}_{m}}\big),
\end{eqnarray*}
so that under $\delta = \delta(T,K,\beta,\alpha)<1$ we obtain for
$q$-a.e. $(s,z)\in[0,T]\times\mathbb{R}$
\begin{eqnarray}\label{recconv2}
\big\|Y^{s,z}\big\|^{2}_{\mathbb{S}^{2}}+\big\|Z^{s,z}\big\|^{2}_{\mathbb{H}^{2}}+
\big\|U^{s,z}\big\|^{2}_{\mathbb{H}^{2}_{m}}\leq C
\mathbb{E}\big[|D_{s,z}\xi|^{2}\big],
\end{eqnarray}
and we arrive at
\begin{eqnarray*}
&&\mathbb{E}\big[\int_{([0,T]\times \mathbb{R})^{2}}|Y^{s,z}(t)|^{2}q(dt,dy)q(ds,dz)\big]<\infty,\nonumber\\
&&\mathbb{E}\big[\int_{([0,T]\times \mathbb{R})^{2}}|Z^{s,z}(t)|^{2}q(dt,dy)q(ds,dz)\big]<\infty,\nonumber\\
&&\mathbb{E}\big[\int_{([0,T]\times
\mathbb{R})^{2}}|U^{s,z}(t,y)|^{2}q(dt,dy)q(ds,dz)\big]<\infty.
\end{eqnarray*}

\noindent Step 4) We show convergence of
$(Y^{n},Z^{n},U^{n})_{n\in\mathbb{N}}$ in
$\mathbb{L}^{1,2}(\mathbb{R})\times\mathbb{L}^{1,2}(\mathbb{R})\times\mathbb{L}^{1,2}(\mathbb{R})$.\\
\noindent From Theorem \ref{thmain} we already know that
$(Y^{n},Z^{n},U^{n})_{n\in\mathbb{N}}$ converges in
$\mathbb{S}^{2}(\mathbb{R})\times\mathbb{H}^{2}(\mathbb{R})\times\mathbb{H}^{2}_{m}(\mathbb{R})$.
We have to prove that the corresponding Malliavin derivatives
converge. The convergence
\begin{eqnarray*}
\lefteqn{\lim_{n\rightarrow\infty}\int_{[0,T]\times\mathbb{R}}\big(\big\|Y^{s,z}-D_{s,z}Y^{n+1}\big\|^{2}_{\mathbb{H}^{2}}}\nonumber\\
&&+\big\|Z^{s,z}-D_{s,z}Z^{n+1}\big\|^{2}_{\mathbb{H}^{2}}+
\big\|U^{s,z}-D_{s,z}U^{n+1}\big\|^{2}_{\mathbb{H}^{2}_{m}}\big)q(ds,dz)=0,
\end{eqnarray*}
for $z=0$ can be proved in the similar way as in the case of a BSDE
without delay driven by a Brownian motion, see for example Theorem
3.3.1 in \cite{I1}.  We only prove the convergence for $z\neq0$.\\
\noindent  Lemma \ref{priori}, applied to the time delayed BSDE
\eqref{eqz}, and \eqref{eqztn1} with \eqref{eqztn2}, yield the
inequality
\begin{eqnarray}\label{fin}
&&\big\|Y^{s,z}-D_{s,z}Y^{n+1}\big\|^{2}_{\mathbb{S}^{2}}+\big\|Z^{s,z}-D_{s,z}Z^{n+1}\big\|^{2}_{\mathbb{H}^{2}}+
\big\|U^{s,z}-D_{s,z}U^{n+1}\big\|^{2}_{\mathbb{H}^{2}_{m}}\nonumber\\
&&\leq (8T+\frac{1}{\beta}) \mathbb{E}\big[\int_{s}^{T}e^{\beta
r}|f^{s,z}(r)-D_{s,z}f^{n}(r)|^{2}dr],
\end{eqnarray}
for $q$-a.e. $(s,z)\in[0,T]\times\mathbb{R}$.\\
\noindent First, by the Lipschitz continuity condition \textbf{(A3)}
for the generator $f$ and the Lipschitz continuity condition
\textbf{(A9)} for the derivative of the function $f$ with respect to
$\omega$ we obtain for $\lambda\otimes m\otimes \lambda$-a.e.
$(s,z,r)\in[0,T]\times(\mathbb{R}-\{0\})\times[0,T]$ the following
two estimates
\begin{eqnarray}\label{derest1}
\lefteqn{|f^{s,z}(r)-D_{s,z}f^{n}(r)|^{2}}\nonumber\\
&\leq&
C\Big(\int_{-T}^{0}|Y^{s,z}(r+v)|^{2}\alpha(dv)+\int_{-T}^{0}|Z^{s,z}(r+v)|^{2}\alpha(dv)\nonumber\\
&&+\int_{-T}^{0}\int_{\mathbb{R}-\{0\}}|U^{s,z}(r+v,y)|^{2}m(dy)\alpha(dv)\nonumber\\
&&\int_{-T}^{0}|D_{s,z}Y^{n}(r+v)|^{2}\alpha(dv)+\int_{-T}^{0}|D_{s,z}Z^{n}(r+v)|^{2}\alpha(dv)\nonumber\\
&&+\int_{-T}^{0}\int_{\mathbb{R}-\{0\}}|D_{s,z}U^{n}(r+v,y)|^{2}m(dy)\alpha(dv)\nonumber\\
&&+\int_{-T}^{0}|Y^{n}(r+v)-Y(r+v)|^{2}\alpha(dv)+\int_{-T}^{0}|Z^{n}(r+v)-Z(r+v)|^{2}\alpha(dv)\nonumber\\
&&+\int_{-T}^{0}\int_{\mathbb{R}-\{0\}}|U^{n}(r+v,y)-U(r+v,y)|^{2}m(dy)\alpha(dv)\Big),
\end{eqnarray}
and for any $\lambda>0$
\begin{eqnarray}\label{derest2}
\lefteqn{|f^{s,z}(r)-D_{s,z}f^{n}(r)|^{2}}\nonumber\\
&\leq&
\big(1+\frac{1}{\lambda}\big)^{2}K\Big(\int_{-T}^{0}|Y^{s,z}(r+v)-D_{s,z}Y^{n}(r+v)|^{2}\alpha(dv)\nonumber\\
&&+\int_{-T}^{0}|Z^{s,z}(r+v)-D_{s,z}Z^{n}(r+v)|^{2}\alpha(dv)\nonumber\\
&&+\int_{-T}^{0}\int_{\mathbb{R}-\{0\}}|U^{s,z}(r+v,y)-D_{s,z}U^{n}(r+v,y)|^{2}m(dy)\alpha(dv)\Big)\nonumber\\
&&+\big(1+\lambda\big)\big(2+\frac{1}{\lambda}\big)K\Big(\int_{-T}^{0}|Y(r+v)-Y^{n}(r+v)|^{2}\alpha(dv)\nonumber\\
&&+\int_{-T}^{0}|Z(r+v)-Z^{n}(r+v)|^{2}\alpha(dv)\nonumber\\
&&+\int_{-T}^{0}\int_{\mathbb{R}-\{0\}}|U(r+v,y)-U^{n}(r+v,y)|^{2}m(dy)\alpha(dv)\Big)/z^{2}.
\end{eqnarray}
\noindent Note that
\begin{eqnarray}\label{inter}
\lefteqn{\int_{[0,T]\times(\mathbb{R}-\{0\})}\big(\big\|Y^{s,z}-D_{s,z}Y^{n+1}\big\|^{2}_{\mathbb{S}^{2}}}\nonumber\\
&&+\big\|Z^{s,z}-D_{s,z}Z^{n+1}\big\|^{2}_{\mathbb{H}^{2}}
+ \big\|U^{s,z}-D_{s,z}U^{n+1}\big\|^{2}_{\mathbb{H}^{2}_{m}}\big)q(ds,dz)\nonumber\\
&=&\lim_{\epsilon\downarrow
0}\int_{0}^{T}\int_{|z|>\epsilon}\big(\big\|Y^{s,z}-D_{s,z}Y^{n+1}\big\|^{2}_{\mathbb{H}^{2}}\nonumber\\
&&+\big\|Z^{s,z}-D_{s,z}Z^{n+1}\big\|^{2}_{\mathbb{H}^{2}}+\big\|U^{s,z}-D_{s,z}U^{n+1}\big\|^{2}_{\mathbb{H}^{2}_{m}}\big)m(dz)ds.
\end{eqnarray}
We prove that this convergence is uniform in $n$.\\
\noindent  Choose $\varepsilon>0$ sufficiently small. By assumption
\textbf{(A7)} we can find $\bar{\epsilon}$ such that
\begin{eqnarray*}
\mathbb{E}\big[\int_{0}^{T}\int_{|z|\leq\bar{\epsilon}}|D_{s,z}\xi|^{2}m(dz)ds\big]<\varepsilon,
\end{eqnarray*}
and
\begin{eqnarray*}
\int_{|z|\leq\bar{\epsilon}}m(dz)<\varepsilon.
\end{eqnarray*}
 Take arbitrary $0<\epsilon_{1}<\epsilon_{2}\leq
\bar{\epsilon}$. By applying the inequality \eqref{fin}, the
estimate \eqref{derest1} and by similar calculations as in
\eqref{contr2} we can derive
\begin{eqnarray}\label{uniform}
\lefteqn{\int_{0}^{T}\int_{\epsilon_{1}<|z|\leq\epsilon_{2}}\big(\big\|Y^{s,z}-D_{s,z}Y^{n+1}\big\|^{2}_{\mathbb{S}^{2}}}\nonumber\\
&&+\big\|Z^{s,z}-D_{s,z}Z^{n+1}\big\|^{2}_{\mathbb{H}^{2}}+\big\|U^{s,z}-D_{s,z}U^{n+1}\big\|^{2}_{\mathbb{H}^{2}_{m}}\big)m(dz)ds\nonumber\\
&\leq& C\int_{0}^{T}\int_{\epsilon_{1}<|z|\leq\epsilon_{2}}
\mathbb{E}\big[\int_{s}^{T}e^{\beta r}|f^{s,z}(r)-D_{s,z}f^{n}(r)|^{2}dr\big]m(dz)ds\nonumber\\
&\leq&
C\Big\{\int_{0}^{T}\int_{\epsilon_{1}<|z|\leq\epsilon_{2}}\big(\big\|Y^{s,z}\big\|^{2}_{\mathbb{S}^{2}}+\big\|Z^{s,z}\big\|^{2}_{\mathbb{H}^{2}}+
\big\|U^{s,z}\big\|^{2}_{\mathbb{H}^{2}_{m}}\big)m(dz)ds\nonumber\\
&&+\int_{0}^{T}\int_{\epsilon_{1}<|z|\leq\epsilon_{2}}\big(\big\|D_{s,z}Y^{n}\big\|^{2}_{\mathbb{S}^{2}}+\big\|D_{s,z}Z^{n}\big\|^{2}_{\mathbb{H}^{2}}+
\big\|D_{s,z}U^{n}\big\|^{2}_{\mathbb{H}^{2}_{m}}\big)m(dz)ds\nonumber\\
&&+\int_{0}^{T}\int_{\epsilon_{1}<|z|\leq\epsilon_{2}}\big(\big\|Y^{n}-Y\big\|^{2}_{\mathbb{S}^{2}}+\big\|Z^{n}-Z\big\|^{2}_{\mathbb{H}^{2}}+
\big\|U^{n}-U\big\|^{2}_{\mathbb{H}^{2}_{m}}\big)m(dz)ds\Big\}.\nonumber\\
\end{eqnarray}
To estimate the first term on the right hand side of
\eqref{uniform}, notice that the inequality \eqref{recconv2} yields
\begin{eqnarray}\label{u1}
\lefteqn{\int_{0}^{T}\int_{\epsilon_{1}<|z|\leq\epsilon_{2}}\big(\big\|Y^{s,z}\big\|^{2}_{\mathbb{S}^{2}}+\big\|Z^{s,z}\big\|^{2}_{\mathbb{H}^{2}}+
\big\|U^{s,z}\big\|^{2}_{\mathbb{H}^{2}_{m}}\big)m(dz)ds}\nonumber\\
&&\leq
C\mathbb{E}\big[\int_{0}^{T}\int_{\epsilon_{1}<|z|\leq\epsilon_{2}}|D_{s,z}\xi|^{2}m(dz)ds\big]<C\varepsilon.
\end{eqnarray}
Recalling $\delta = \delta(T,K,\beta,\alpha)<1$ and applying the
inequality \eqref{recconv1} we estimate the second term in
\eqref{uniform} by
\begin{eqnarray}\label{u2}
\lefteqn{\int_{0}^{T}\int_{\epsilon_{1}<|z|\leq\epsilon_{2}}\big(\big\|D_{s,z}Y^{n}\big\|^{2}_{\mathbb{S}^{2}}+\big\|D_{s,z}Z^{n}\big\|^{2}_{\mathbb{H}^{2}}+
\big\|D_{s,z}U^{n}\big\|^{2}_{\mathbb{H}^{2}_{m}}\big)m(dz)ds}\nonumber\\
&\leq&9e^{\beta T}\mathbb{E}\big[\int_{0}^{T}\int_{\epsilon_{1}<|z|\leq\epsilon_{2}}|D_{s,z}\xi|^{2}m(dz)ds\big]\nonumber\\
&&+\delta\int_{0}^{T}\int_{\epsilon_{1}<|z|\leq\epsilon_{2}}\big(\big\|D_{s,z}Y^{n-1}\big\|^{2}_{\mathbb{S}^{2}}+\big\|D_{s,z}Z^{n-1}\big\|^{2}_{\mathbb{H}^{2}}+
\big\|D_{s,z}U^{n-1}\big\|^{2}_{\mathbb{H}^{2}_{m}}\big)m(dz)ds\nonumber\\
&<&\frac{9e^{\beta T}\varepsilon}{1-\delta}\nonumber\\
&&+\delta^{n}\int_{0}^{T}\int_{\epsilon_{1}<|z|\leq\epsilon_{2}}\big(\big\|D_{s,z}Y^{0}\big\|^{2}_{\mathbb{S}^{2}}+\big\|D_{s,z}Z^{0}\big\|^{2}_{\mathbb{H}^{2}}+
\big\|D_{s,z}U^{0}\big\|^{2}_{\mathbb{H}^{2}_{m}}\big)m(dz)ds.
\end{eqnarray}
The estimate of the third term follows from the contraction
inequality \eqref{contr}
\begin{eqnarray}\label{u3}
\lefteqn{\int_{0}^{T}\int_{\epsilon_{1}<|z|\leq\epsilon_{2}}\big(\big\|Y^{n}-Y\big\|^{2}_{\mathbb{S}^{2}}+\big\|Z^{n}-Z\big\|^{2}_{\mathbb{H}^{2}}+
\big\|U^{n}-U\big\|^{2}_{\mathbb{H}^{2}_{m}}\big)m(dz)ds}\nonumber\\
&&\quad \ \quad \leq
\delta^{n}T\big(\big\|Y^{0}-Y\big\|^{2}_{\mathbb{S}^{2}}+\big\|Z^{0}-Z\big\|^{2}_{\mathbb{H}^{2}}+
\big\|U^{0}-U\big\|^{2}_{\mathbb{H}^{2}_{m}}\big)\varepsilon.
\end{eqnarray}
Choosing $Y^{0}=Z^{0}=U^{0}=0$ and combining \eqref{u1}-\eqref{u3} gives the uniform convergence of \eqref{inter}.\\
\noindent Next, by applying the inequality \eqref{fin}, the estimate
\eqref{derest2} and similar calculations as in \eqref{contr2} and
\eqref{contr} we can derive
\begin{eqnarray*}
\lefteqn{\int_{0}^{T}\int_{|z|>\epsilon}\big(\big\|Y^{s,z}-D_{s,z}Y^{n+1}\big\|^{2}_{\mathbb{S}^{2}}}\nonumber\\
&&+\big\|Z^{s,z}-D_{s,z}Z^{n+1}\big\|^{2}_{\mathbb{H}^{2}}+\big\|U^{s,z}-D_{s,z}U^{n+1}\big\|^{2}_{\mathbb{H}^{2}_{m}}\big)m(dz)ds\nonumber\\
&\leq& \big(8T+\frac{1}{\beta}\big)\int_{0}^{T}\int_{|z|>\epsilon}
\mathbb{E}\big[\int_{s}^{T}e^{\beta r}|f^{s,z}(r)-D_{s,z}f^{n}(r)|^{2}dr]m(dz)ds\nonumber\\
&\leq&
\delta\Big\{\big(1+\frac{1}{\lambda}\big)^{2}\int_{0}^{T}\int_{|z|>\epsilon}\big(\big\|Y^{s,z}-D_{t,z}Y^{n}\big\|^{2}_{\mathbb{S}^{2}}\nonumber\\
&&+\big\|Z^{s,z}-D_{s,z}Z^{n}\big\|^{2}_{\mathbb{H}^{2}}+
\big\|U^{s,z}-D_{s,z}U^{n}\big\|^{2}_{\mathbb{H}^{2}_{m}}\big)m(dz)ds\nonumber\\
&&+\big(1+\lambda\big)\big(2+\frac{1}{\lambda}\big)\big(\big\|Y^{n}-Y\big\|^{2}_{\mathbb{S}^{2}}+\big\|Z^{n}-Z\big\|^{2}_{\mathbb{H}^{2}}+
\big\|U^{n}-U\big\|^{2}_{\mathbb{H}^{2}_{m}}\big)\int_{|z|>\epsilon}\nu(dz)\Big\},
\end{eqnarray*}
and we choose $\lambda$ sufficiently large such that
$\tilde{\delta}:=\delta\big(1+\frac{1}{\lambda}\big)^{2}<1$.\\
\noindent Due to the convergence of
$(Y^{n},Z^{n},U^{n})_{n\in\mathbb{N}}$, for a sufficiently small
$\varepsilon>0$ we can find $N$ sufficiently large such that for all
$n\geq N$
\begin{eqnarray*}
\big(1+\lambda\big)\big(2+\frac{1}{\lambda}\big)\big(\big\|Y^{n}-Y\big\|^{2}_{\mathbb{S}^{2}}+\big\|Z^{n}-Z\big\|^{2}_{\mathbb{H}^{2}}+
\big\|U^{n}-U\big\|^{2}_{\mathbb{H}^{2}_{m}}\big)\int_{|z|>\epsilon}\nu(dz)<\varepsilon.
\end{eqnarray*}
We derive the recursion for $n\geq N$
\begin{eqnarray*}
\lefteqn{\int_{0}^{T}\int_{|z|>\epsilon}\big(\big\|Y^{s,z}-D_{s,z}Y^{n+1}\big\|^{2}_{\mathbb{S}^{2}}}\nonumber\\
&&+\big\|Z^{s,z}-D_{s,z}Z^{n+1}\big\|^{2}_{\mathbb{H}^{2}}+\big\|U^{s,z}-D_{s,z}U^{n+1}\big\|^{2}_{\mathbb{H}^{2}_{m}}\big)m(dz)ds\nonumber\\
&<&\tilde{\delta}\big\{\int_{0}^{T}\int_{|z|>\epsilon}\big(\big\|Y^{s,z}-D_{s,z}Y^{n}\big\|^{2}_{\mathbb{S}^{2}}\nonumber\\
&&+\big\|Z^{s,z}-D_{s,z}Z^{n}\big\|^{2}_{\mathbb{H}^{2}}+\big\|U^{s,z}-D_{s,z}U^{n}\big\|^{2}_{\mathbb{H}^{2}_{m}}\big)m(dz)ds\}+\delta\varepsilon\nonumber\\
&<&\tilde{\delta}^{n-N}\int_{0}^{T}\int_{|z|>\epsilon}\big(\big\|Y^{s,z}-D_{s,z}Y^{N}\big\|^{2}_{\mathbb{S}^{2}}\nonumber\\
&&+\big\|Z^{s,z}-D_{s,z}Z^{N}\big\|^{2}_{\mathbb{H}^{2}}+\big\|U^{s,z}-D_{s,z}U^{N}\big\|^{2}_{\mathbb{H}^{2}_{m}}\big)m(dz)ds+\frac{\delta\varepsilon}{1-\tilde{\delta}},
\end{eqnarray*}
and finally we conclude that
\begin{eqnarray*}
\lefteqn{\lim_{n\rightarrow\infty}\int_{0}^{T}\int_{|z|>\epsilon}\big(\big\|Y^{s,z}-D_{s,z}Y^{n+1}\big\|^{2}_{\mathbb{S}^{2}}}\nonumber\\
&&+\big\|Z^{s,z}-D_{s,z}Z^{n+1}\big\|^{2}_{\mathbb{H}^{2}}+\big\|U^{s,z}-D_{s,z}U^{n+1}\big\|^{2}_{\mathbb{H}^{2}_{m}}\big)m(dz)ds=0.
\end{eqnarray*}
The equation
\begin{eqnarray*}
\lefteqn{\lim_{n\rightarrow\infty}\int_{[0,T]\times(\mathbb{R}-\{0\})}\big(\big\|Y^{s,z}-D_{s,z}Y^{n+1}\big\|^{2}_{\mathbb{S}^{2}}}\nonumber\\
&&+\big\|Z^{s,z}-D_{s,z}Z^{n+1}\big\|^{2}_{\mathbb{H}^{2}}+\big\|U^{s,z}-D_{s,z}U^{n+1}\big\|^{2}_{\mathbb{H}^{2}_{m}}\big)q(ds,dz)=0
\end{eqnarray*}
now follows by interchanging the limits in $n$ and $\varepsilon$ in \eqref{inter}.\\

\noindent Step 4) Since the space $\mathbb{L}^{1,2}(\mathbb{R})$ is
a Hilbert space and the Malliavin derivative is a closed operator,
see Theorem 12.6 in \cite{DIN}, the claim that
$(Y,Z,U)\in\mathbb{L}^{1,2}(\mathbb{R})\times\mathbb{L}^{1,2}(\mathbb{R})\times\mathbb{L}^{1,2}(\mathbb{R})$
and $(Y^{s,z}(t),Z^{s,z}(t),U^{s,z}(t,y))_{0\leq s,t\leq
T,(y,z)\in(\mathbb{R}-\{0\})\mathbb{R}}$ is a version of the
derivative $(D_{s,z}Y(t),D_{s,z}Z(t),D_{s,z}U(t,y))_{0\leq s,t\leq
T,(y,z)\in(\mathbb{R}-\{0\})\mathbb{R}}$ follows, and finishes the
proof. \cbdu

The following Corollary shows that the interpretation of the
solution component $(Z,U)$ in terms of the Malliavin trace of $Y$
still holds for BSDE with time delayed generators.

\begin{cor}
Under the assumptions of Theorem \ref{thmmalliavin}, we have
\begin{center}
$\big((D_{t,0}Y)^{\mathcal{P}}(t)\big)_{0\leq t \leq T}$ is a version of $\big(Z(t)\big)_{0\leq t\leq T}$,\\
$\big((D_{t,z}Y)^{\mathcal{P}}(t)\big)_{0\leq t\leq T,
z\in(\mathbb{R}-\{0\})}$ is a version of $\big(U(t,z)\big)_{0\leq
t\leq T, z\in(\mathbb{R}-\{0\})}$,
\end{center}
where $(\cdot)^{\mathcal{P}}$ denotes the predictable projection of
a process.
\end{cor}
\Proof The solution of \eqref{bsdem} satisfies
\begin{eqnarray}\label{pp}
\lefteqn{Y(s)=Y(0)}\nonumber\\
&&-\int_{0}^{s}f\Big(r,\int_{-T}^{0}Y(r+v)\alpha(dv),\nonumber\\
&&\quad \ \quad \int_{-T}^{0}Z(r+v)\alpha(dv), \int_{-T}^{0}\int_{\mathbb{R}-\{0\}}U(r+v,y)m(dy)dv\Big)dr\nonumber\\
&&+\int_{0}^{s}Z(r)dW(r)+\int_{0}^{s}\int_{\mathbb{R}-\{0\}}U(r,y)\tilde{M}(dr,dy),\quad
0\leq s\leq T.
\end{eqnarray}
By differentiating \eqref{pp} we obtain according to Lemma
\ref{derl} for $q$-a.e. $(u,z)\in[0,T]\times\mathbb{R}$
\begin{eqnarray*}
D_{u,0}Y(s)&=&Z(u)-\int_{u}^{s}D_{u,0}f(r)dr
+\int_{u}^{s}D_{u,0}Z(r)dW(r)\nonumber\\
&&+\int_{u}^{s}\int_{\mathbb{R}-\{0\}}D_{u,0}U(r,y)\tilde{M}(dr,dy),\quad
0\leq u\leq s\leq T,
\end{eqnarray*}
and for $z\neq0$
\begin{eqnarray*}
D_{u,z}Y(s)&=&U(u,z)-\int_{u}^{s}D_{u,z}f(r)dr+\int_{u}^{s}D_{u,z}Z(r)dW(r)\nonumber\\
&&+\int_{u}^{s}\int_{\mathbb{R}-\{0\}}D_{u,z}U(r,y)\tilde{M}(dr,dy),\quad
0\leq u\leq s\leq T,
\end{eqnarray*}
where the derivative operators $D_{u,z}$ are defined by\eqref{op0}
and \eqref{opz}. Since the mappings $s\mapsto\int_{u}^{s}
D_{u,z}f(r)dr$,$s\mapsto\int_{u}^{s}D_{u,z}Z(r)dW(r)$ are
$\mathbb{P}$-a.s. continuous and the mapping $s\mapsto
\int_{u}^{s}\int_{\mathbb{R}-\{0\}}D_{u,z}U(r,y)\tilde{M}(dr,dy)$ is
$\mathbb{P}$-a.s. c\`{a}dl\`{a}g (see Theorems 4.2.12 and 4.2.14 in
\cite{A}), taking the limit $s\downarrow u$ yields
\begin{eqnarray*}
D_{u,0}Y(u)&=&Z(u),\,\,\mbox{for}\,\, \lambda-a.e. u\in[0,T],
\mathbb{P}-a.s.,\nonumber\\
D_{u,z}Y(u)&=&U(u,z)\,\,\mbox{for}\,\, \lambda\otimes m-a.e.
(u,z)\in[0,T]\times(\mathbb{R}-\{0\}), \mathbb{P}-a.s..
\end{eqnarray*}
As $Y\in\mathbb{S}^{2}(\mathbb{R})$ has $\mathbb{P}-a.s.$
c\`{a}dl\`{a}g $\mathbb{F}$-adapted trajectories, for $0\leq u\leq
T$ we have the representation
\begin{eqnarray*}
Y(u)=\sum_{n=0}^{\infty}I_{n}(g_{n}((u,0),.)=\sum_{n=0}^{\infty}I_{n}(g_{n}((u,0),.)\mathbf{1}^{\otimes
n}_{[0,u]}(.)),\quad g_{n}\in L^{2}_{T,q,n+1},n\geq 0,
\end{eqnarray*}
with c\`{a}dl\`{a}g mappings $u\mapsto g_{n}((u,0),.)$. By
Definition 3.1.2 of the Malliavin derivative we arrive at
\begin{eqnarray*}
D_{u,z}Y(u)=\sum_{n=0}^{\infty}nI_{n-1}(g_{n}((u,0),(u,z),.)\mathbf{1}^{\otimes
n}_{[0,u]}((u,z),.)),\,\,\mbox{for}\,\, q-a.e.
(u,z)\in[0,T]\times\mathbb{R}.
\end{eqnarray*}
For $\delta_{\{0\}}\times m$-a.e. $z\in\mathbb{R}$, we conclude that
the mapping $(u,\omega)\mapsto D_{u,z}Y(u)(\omega)$ is
$\mathbb{F}$-adapted and measurable and has a progressively
measurable (optional) modification. Moreover, notice that the
optional process $u\mapsto D_{u,z}Y(u)$ and its unique predictable
projection $u\mapsto (D_{u,z}Y)^{\mathcal{P}}(u)$ are modifications
of each other, see Theorem 5.5 in \cite{Y}. Finally, we remark that
there exists a $\mathcal{P}\times\mathcal{B}(\mathbb{R})$ measurable
version of $(\omega,u,z)\mapsto(D_{u,z}Y)^{\mathcal{P}}(u)(\omega)$,
see Lemma 2.2 in \cite{AI}. This completes the proof. \cbdu
\\
\
\\
\noindent \textbf{Acknowledgements:} This paper was written while
the first author was staying at Humboldt University at Berlin.
{\L}ukasz Delong acknowledges the financial support from the AMaMeF
programme.

\end{document}